\documentclass[12pt]{amsart}
\usepackage{amsthm,amsfonts,amssymb,amscd}

\usepackage[all]{xy}

\newcommand\D{{\bf{D}}}
\newcommand\F{{\boldsymbol{\mathcal F}}}
\newcommand\bH{{\bf{H}}}
\newcommand\J{{\bf{J}}}
\newcommand\X{{\bf X}}
\newcommand\Y{{\bf{Y}}}

\newcommand\aalpha{{\boldsymbol{\alpha}}}
\newcommand\ppi{{\boldsymbol{\pi}}}
\newcommand\SSigma{{\boldsymbol{\Sigma}}}

\newcommand{\CC}{{\mathbb C}}

\newcommand{\PP}{{\mathbb P}}

\newcommand{\BBB}{{\mathcal B}}
\newcommand{\DDD}{{\mathcal D}}
\newcommand{\OOO}{{\mathcal O}}

\newcommand{\III}{{\mathcal I}}
\newcommand{\JJJ}{{\mathcal J}}

\newcommand{\EEE}{{\mathcal E}}
\newcommand{\HHH}{{\mathcal H}}
\newcommand{\LLL}{{\mathcal L}}
\newcommand{\FFF}{{\mathcal F}}
\newcommand{\CCC}{{\mathcal C}}

\newcommand{\NNN}{{\mathcal N}}

\newcommand{\SSS}{{\mathcal S}}

\newcommand{\EXT}{{\mathcal Ext}}

\newcommand\uC{{\underline{C}}}

\newcommand\codim{\text{codim}}

\newcommand{\const}{\operatorname{const}\nolimits}

\newcommand\Ext{\operatorname{Ext}\nolimits}

\newcommand{\Hilb}{\operatorname{Hilb}\nolimits}

\newcommand\im{\operatorname{Im}\nolimits}
\newcommand\Ker{\operatorname{Ker}\nolimits}

\newcommand\pr{\text{pr}}
\newcommand{\Proj}{\operatorname{Proj}\nolimits}
\newcommand\red{\text{red}}
\newcommand\res{\text{res}}

\newcommand{\Sing}{\operatorname{Sing}\nolimits}
\newcommand\Span{\text{Span}}

\newcommand\Sym{\operatorname{Sym}\nolimits}

\newcommand{\ka}{{\mathcal A}}

\newcommand{\kc}{{\mathcal C}}
\newcommand{\kd}{{\mathcal D}}

\newcommand{\kf}{{\mathcal F}}

\newcommand{\ki}{{\mathcal I}}

\newcommand{\kl}{{\mathcal L}}

\newcommand{\ko}{{\mathcal O}}

\newcommand{\ks}{{\mathcal S}}

\newcommand\elm{{\operatorname{elm}\nolimits}}

\newtheorem{theorem}{Theorem}[section]
\newtheorem{proposition}[theorem]{Proposition}

\newtheorem{lemma}[theorem]{Lemma}
\newtheorem{corollary}[theorem]{Corollary}
\newtheorem{sub}[subsection]{}

\theoremstyle{definition}
\newtheorem{definition}[theorem]{Definition}

\newtheorem{proposition-definition}[theorem]{Proposition-Definition}

\theoremstyle{remark}

\newcommand{\into}{\hookrightarrow}
\newcommand{\tto}[1]{\xrightarrow{#1}}

\newcommand\lra{{\longrightarrow}}

\newcommand\rar{\rightarrow}

\newlength{\rrrr}
\newcommand{\isom}[1]{{\settowidth{\rrrr}{$\scriptstyle{x#1x}$}
\xrightarrow{\makebox[\rrrr]{$\scriptstyle{#1}$}}
\hspace{-0.5\rrrr }\hspace{-1.1 em}
\raisebox{- 0.5 ex}{$\sim$}\hspace{0.7\rrrr }
}}

\renewcommand\square{\frame{\phantom{{\large x}}}}
\newcommand\empt{\varnothing}
\renewcommand\emptyset{\varnothing}

\begin{document}

\title{A parametrization of the theta divisor of the quartic double solid\\
}

\author[Markushevich]{\;D.G~Markushevich}

\address{D. M.: Math\'ematiques - b\^{a}t. M2, Universit\'e Lille 1,
F-59655 Villeneuve d'Ascq Cedex, France}
\email{markushe@agat.univ-lille1.fr}

\author[Tikhomirov]{\;A.S.~Tikhomirov}

\address{
  Department of Mathematics\\
  State Pedagogical University\\
  Respublikanskaya Str. 108
\newline 150 000 Yaroslavl, Russia}
\email{tikhomir@yaroslavl.ru}

\thanks{Partially supported by the grant
INTAS-OPEN-2000-269}

\begin{abstract}
Let $M_X(2;0,3)$ be the moduli space of rank-2 stable
vector bundles with Chern classes $c_1=0, c_2=3$ on the Fano threefold $X$,
the double solid $\mathbb P^3$ of index two. We prove that the vector bundles
obtained by Serre's construction from smooth elliptic quintic curves on $X$
form an open part of an irreducible component $M$ of $M_X(2;0,3)$
and that the Abel-Jacobi map $\phi:M\to J(X)$ into the intermediate Jacobian
$J(X)$ of $X$ defined by the second Chern class
is generically finite of degree 84 onto a translate $\Theta+\const$
of the theta-divisor. We also prove that the family of elliptic quintics
on a general $X$ is irreducible and of dimension 10.
\end{abstract}


\maketitle
\thispagestyle{empty}


\setcounter{section}{-1}

\section{Introduction}

After the famous paper of Clemens and Griffiths \cite{CG},
in which they determined the intermediate Jacobian
and its theta divisor for the cubic threefold,
the next Fano threefold to study was naturally the quartic double solid,
that is the double cover \mbox{$\pi:X_2\xrightarrow{2:1}\PP^3$}
ramified in a quartic surface $W\subset\PP^3$.
In 1981, Welters \cite{We} found a family of curves
parametrizing the intermediate Jacobian $J(X_2)$,
the septics of genus 4, and this is the simplest known family
of generically irreducible curves with this property. He developed
also techniques that allow one to decide whether the
Abel-Jacobi images of certain families of curves lie in
a translate $\Theta+\const$ of the theta divisor, but failed to find one
parametrizing the whole of $\Theta+\const$. Such a family
was found in 1986 (with a use of the techniques of Welters)
by Tikhomirov \cite{T-2}: these are
the sextics of genus 3, or the so called Reye sextics.
In 1991, Clemens \cite{C-2} reproved Tikhomirov's result
by degenerating the double solid into a pair
of $\PP^3$'s meeting each other transversely along a quadric,
a technique he developed earlier in \cite{C-1}.
Several other papers have contributed, since 1981, to the
study of quartic double solids: \cite{T-1}, \cite{SV},
\cite{De}, \cite{V}, \cite{Iz}.

In 2001, Tikhomirov \cite{T-3} found a new and simpler parametrization
of $\Theta$,
the one by the family of elliptic quintics $\CCC_5^1(X)$.
Its advantage is that the elliptic quintics
define stable vector bundles on $X$ via Serre's construction,
and the Abel-Jacobi map factors through the
moduli space of vector bundles.

The study of moduli spaces of stable vector bundles
on Fano threefolds of indices 1 and 2 is quite a
recent topic. Recall that the {\em index} of a Fano
threefold $X$ is the maximal
integer $\nu$ dividing $K_X$ in the Picard group of $X$.
The cubic $X_3$ and the quartic double solid $X_2$ are Fano
of index 2. The nature of the results obtained
depends strongly on the index.
In the index-2 case, the answers are given in terms
of the Abel--Jacobi map defined by the
second Chern class $c_2$ modulo rational equivalence.
In the index-1 case, the Abel--Jacobi map does not bring
much new information about the moduli spaces.

The results known so far include
the desciption of one component of a moduli space $M_{X}(2;c_1,c_2)$
of rank-2 vector bundles with $c_1=0$ or 1 and small $c_2$
on each one of the following five Fano threefolds:
the cubic $X_3$ \cite{MT-1}, \cite{IM-1}, \cite{Druel}, \cite{B},
the quartic \cite{IM-2}, the Fano threefolds $X_{2g-2}$ of the main series
of genus 7 \cite{IM-3} and 9 \cite{IR},
and the double solid $X_2$, treated in the present paper and in \cite{T-3}.
It is noteworthy that in most cases the interesting
components of moduli are those parametrizing the vector
bundles obtained by Serre's construction from elliptic quintics.
For the cubic $X_3$, the relevant moduli space
is $M_{X_3}(2;0,2)$, and the Abel-Jacobi map identifies it with
an open set in the intermediate Jacobian $J(X_3)$.
For the double solid,
the corresponding moduli space
$M_{X}(2;0,3)$ is reducible: it has at least two components, one of which,
say $M$, is associated to
elliptic quintics, and the other to disjoint unions
of two elliptic curves in $X$ whose images in $\PP^3$ are
a plane cubic and a line. Leaving apart the question whether
there are any other components, we restrict ourselves
to the study of $M$.

The main result of this paper states that the Abel-Jacobi
map $\Phi$ sends the family of elliptic quintics $\CCC_5^1(X)$
onto an open subset of $\Theta+\const$,
and that the Serre construction defines
a factorization $\CCC_5^1(X)\xrightarrow{\mbox{\scriptsize Serre}}
M\xrightarrow{g}\Theta+\const$ with $g$ generically finite
of degree 84 (Theorems \ref{Theta} and \ref{Stein}).
We also prove the irreducibility of $\CCC_5^1(X)$ (Theorem \ref{quintics}),
so that $M$ is really an irreducible component
of the moduli space, rather than a union of components.
The facts that the image of $\Phi$ is dense in $\Theta+\const$
and that $g$ is generically finite were stated in \cite{T-3}, so
the main goal of our paper is to determine the degree of $g$.
Nonetheless, we provide a complete proof of these facts too,
for their proof was only sketched in loc. cit.

We use Welters' criterion (Proposition \ref{Welters}) which gives sufficient
conditions for a certain  subfamily $\{ C_t\}_{t\in T_l}$
of a given family $\{ C_t\}_{t\in T}$
of curves of degree $k(k+1)-3$ to be mapped
by the Abel-Jacobi map into a translate of
the theta divisor. Here $l$ is a line
in $X$ and
$$
T_l=\{t\in T\mid C_t\cap l=\empt~\mbox{ and}~C_t\cup
l~\mbox{ lies~on~a~surface~from}~|\OOO_X(k)|\} ,
$$
where $\OOO_X(k):=\pi^*\OOO_{\PP^3}(k)$. We construct
such a family $T$ with $k=4$ as a partial compactification of the family
of reducible curves
$C=C_1\cup C_2$, where $C_1$ is a genus-4 septic, $C_2$
the pullback of a genus 2 quintic $\uC_2\subset\PP^3$,
and $C_1,C_2$ meet each other quasi-transversely in 13
points. We show that $T_l$ contains a subfamily of curves of the form
$C'\cup C''\cup C_2$ such that $C'\in\CCC_5^1(X)$
and $C''=\pi^{-1}(m)$, where $m\subset\PP^3$ is a trisecant
of $\pi (C')$, whose Abel-Jacobi images differ from those
of $C'$ by a constant translation. Hence
$\Phi_{\CCC_5^1(X)} (\CCC_5^1(X))\subset
\Phi_{T_l} (T_l)+\const \subset\Theta+\const$.
The opposite inclusion follows from
the computation of the differential of the Abel-Jacobi map.
The generic fiber of $\Phi_{\CCC_5^1(X)}$ is a union of
$\nu$ copies of $\PP^1$, which can be interpreted either as pencils of
elliptic quintics $\{ C_\lambda\}_{\lambda\in\PP^1}$
in the K3 surfaces $S(C_\lambda )$ from the linear
system $|\OOO_X(2)|$, or as pencils of the zero loci
of sections $s\in H^0(X,\EEE )$, where $\EEE$ is the vector bundle
obtained by Serre's construction from $C_\lambda$.

To compute the number $\nu$ of copies of $\PP^1$ in the generic
fiber of $\Phi_{\CCC_5^1(X)}$, we use the degeneration $\{ X_t\}_{t\in\Delta}$
of $X$ into the union $X_0=\PP^{3'}\cup\PP^{3''}$ of two copies of $\PP^3$
meeting each other in a smooth quadric $Q=\{ G=0\}$, where $\Delta$ is a disc
in $\CC$. The fiber $X_t$
is defined as the double cover $\pi_t:X_t\xrightarrow{2:1}\PP^3$
ramified in $W_t=\{ tF+G^2=0\}$, where $W=\{ F=0\}$ is the branch
locus of $X\lra\PP^3$. The Hodge theory and the Neron model
of the relative intermediate Jacobian of such a degeneration
were studied by Clemens in \cite{C-1}.
He made also the following important observation: let $C_0\in X_0$
be an irreducible curve of degree $d$ which is the limit of a family of
curves $C_s\subset X_t$, $t=s^e$. Assume that $C_s$ for $s\neq 0$
is projected birationally onto its image $\uC_s\subset\PP^3$; this means
that $\uC_s$ is tangent to $W_t$ in $2d$ points. Then $C_0$ is
a $2d$-secant to the octic curve $B=W\cap Q$. Thus the limits of
elliptic quintics from $\CCC_5^1(X_t)$ are elliptic quintics
in $\PP^{3 '}$ or $\PP^{3''}$ which are 10-secant to $B$. The work of
Clemens allows us to extend the Abel-Jacobi map to such
curves in $X_0$ and the question on the number of copies
of $\PP^1$ in the generic fiber of the Abel-Jacobi map $\Phi_0$ at $t=0$ is
reduced to that on the number of elliptic quintics in $\PP^3$
which are 10-secant to a given octic $B$.

By a result of Getzler \cite{Ge}, there are 42
elliptic quintics passing through
10 generic points in $\PP^3$. We show that the 10-uples of points
lying in the smooth complete intersections
of a quartic and a quadric in $\PP^3$ form a divisor in the main component
$\HHH$ of $\Hilb^{10}(\PP^3)$ and that this divisor
is not contained in the branch
locus of the 42-sheeted covering defined by the elliptic quintics,
hence there are $2\cdot 42=84$ copies of $\PP^1$ in the generic
fiber of $\Phi_0$ (42 copies in each one of the
two $\PP^3$'s constituting $X_0$). We show also that the elliptic
quintics, 10-secant to $B$, are acquired with multiplicity 1
in the relative Hilbert scheme of the family $\{ X_t\}_{t\in\Delta}$,
hence the number of copies of $\PP^1$ in the fiber of $\Phi_t$
is invariant under the deformation and $\nu =84$.

To this end, we show that the family of 10-secant elliptic quintics
to a given generic $B$ is irreducible and construct a family
$\{ C_t\}_{t\in\Delta}$ of cycles of 5 lines,
or pentagons in $X_t$, such that
$C_0$ is 10-secant to $B$ and is strongly smoothable into
an elliptic quintic
in $X_0$, 10-secant to $B$.
This implies that a generic 10-secant elliptic quintic in $X_0$
admits a local cross-section passing through it in the relative
Hilbert scheme of $\{ X_t\}_{t\in\Delta}$, hence the fiber over
$t=0$ is acquired with multiplicity 1.

We now describe the contents of the paper by sections.

The main result of Section 1 is the proof of the irreducibility of the
family of elliptic quintics $\CCC_5^1(X)$.
We start by proving the irreducibility of the families $\CCC_2^0(X)$ of
conics and $\CCC_3^0(X)$ of twisted cubics.
The proofs are based on the splitting of the curves under
consideration into two components of smaller degrees.
As concerns elliptic quintics,
we split them into a twisted cubic and an "elliptic
conic", that is, $\pi^{-1}(m)$ for a line $m\subset\PP^3$.
As soon as the irreducibility of the family of split curves
is established, we get a distinguished component of $C_5^1(X)$, namely
the one containing the smoothings of the split curves. Then we use
the irreducibility of the monodromy action on the set
of all the components of $C_5^1(X)$ to deduce that the number of components
is 1.

In Section 2 we prove the irreducibility of the family
$C_7^4(X)$ of genus-4 septics in $X$, which dominate
$J(X)$ according to Welters, and construct the above mentioned family
$T$ of curves of degree $k(k+1)-3=17$ for $k=4$,
to which the criterion of Welters is applied.
We conclude by the proof of Theorem \ref{Theta},
stating that the Abel-Jacobi image of the elliptic quintics
is a translate of $\Theta$.

In Section 3 we explain the Clemens' degeneration technique
and describe the limiting Abel-Jacobi map $\Phi_0$ on the degenerate
double solid, which consists of two copies of $\PP^3$.
We also show that the central fiber in the domain
of the Abel-Jacobi map is reduced by studying the
deformations of 10-secant pentagons. This implies that
the number of connected components in the generic
fiber of $\Phi_t$ is the same as for $\Phi_0$.

In Section 4 we prove that the number of elliptic quintics
in $\PP^3$ passing through 10 generic points
remains unchanged if we specialize the 10-uple of points
to a 10-uple lying on a generic complete intersection
of a quadric and a quartic. This implies that
the number of connected components in the generic
fiber of $\Phi_0$ is 84.

In Section 5 we introduce the Serre construction and
reformulate the main theorem in terms of moduli of
vector bundles.

\bigskip

{\sc Acknowledgements.} The authors thank the referee,
whose remarks helped them to improve the presentation.

\section{Conics, twisted cubics and elliptic quintics on $X$}
\label{cocuqu}

\begin{sub}\label{nota-irred}\rm {\bf Notation.}
\begin{itemize}
\item $X$, a general quartic double solid, $X\subset
Y={\Proj}(\ko_{{\PP}^3}\oplus\ko_{{\PP}^3}(2))
$,
where $\kl=\ko_{Y/{\PP}^3}(1)$ (see \cite{T-1} or \cite{We}).
\item $\rho :Y\lra\PP^3$, $\pi =\rho|_X:X\lra\PP^3$,
the natural maps.
\item $\kf=\kf(X)$, the Fano surface of $X$,
that is the base of the family of lines on~$X$.
\item $\mathcal{C}^0_2=\mathcal{C}^0_2(X)\subset \Hilb(X)$,
the base of the family of smooth conics in $X$.
\item $\kc^1_2(X):=
\{C\in\Hilb(X)~|~C=\pi^{-1}(l),\ l\in G(1,3)\}$; we call the curves
from $\kc^1_2(X)$ elliptic conics on $X$.
\item $\kc_3^0(\PP^3)$, $\mathcal{C}^0_3=\mathcal{C}^0_3(X)$,
resp. $\mathcal{C}^0_3(Y)$, the base of the family of
twisted cubics in $\PP^3$, $X$, resp. $Y$.
\item $\kc_5^1(\PP^3)$,
$\mathcal{C}^1_5=\mathcal{C}^1_5(X)$), resp. $\mathcal{C}^1_5(y)$,
the base of the family of
elliptic quintics in $\PP^3$, $X$, resp. $Y$.
\end{itemize}
A line (resp. conic, twisted cubic, elliptic quintic)
in $Y$ is, by definition, a curve $C\subset Y$
such that $\rho |_C$ is an isomorphism onto a line
(resp. conic, rational
twisted cubic, elliptic quintic) $C_0=\rho (C)\subset\PP^3$.
In general, we use the notation $\CCC_d^g(X)$ for the base of
a family of curves of degree $d$ and genus $g$ on $X$, the degree on $X$
being defined with respect to the ample sheaf $\OOO_X(1):=
\pi^*\OOO_{\PP^3}(1)$.

\end{sub}

The main result of this section is the irreducibility of
the family of elliptic quintics on $X$ (Theorem \ref{quintics}).
First we need to prove
that the families $\mathcal{C}^0_2(X)$
and $\kc^0_3(X)$ of conics and, respectively,
twisted cubics on $X$ are irreducible.

\begin{sub}\label{Conics} \rm {\bf Conics.}
First recall the well known facts about lines on a general X
(see, e.g., \cite{T-1},
\cite{We}).
\begin{lemma}\label{Lines}
Let $X$ be a general quartic double solid. Then the following assertions hold:

(i) The Fano surface $\kf=\kf(X)$ of lines on $X$ is an irreducible surface
having at most isolated quadratic points as singularities, and
the incidence divisor $D_\kf=\{(l_1,l_2)\in\kf\times\kf~|~l_1\cap l_2\ne\empt\}$
on $\kf\times\kf$ is irreducible.

(ii) For a general line $l\in \mathcal{F}(X)$, the normal bundle
$N_{l/X}$ is isomorphic to $2\mathcal{O}_{{\PP}^1}$.
\end{lemma}

Next we will prove the following statement:

\begin{lemma}\label{conics}
Let $X$ be a general quartic double solid. Then:

(i) $\mathcal{C}^0_2(X)$ is irreducible of dimension $4$.

(ii) For a general conic $C\in\mathcal{C}^0_2(X)$, the normal bundle
$N_{C/X}$ is isomorphic either to $2\mathcal{O}_{{\PP}^1}(1)$ or to
$\mathcal{O}_{{\PP}^1}\oplus\mathcal{O}_{{\PP}^1}(2)$.
\end{lemma}

\begin{proof}
(i) Let $\check W\subset\check\PP^3$ be the dual surface of the quartic $W$
and consider the morphism
$v:\mathcal{C}^0_2(X)\to\check\PP^3:\ C\mapsto\Span(\pi(C))$.
Then
$v^{-1}(Y)$
for a plane
$Y\in\check\PP^3$
consists of all conics on the surface
$S_Y=\pi^{-1}(Y)$.
If
$Y\in U:=\check\PP^3\smallsetminus\check W$,
then $S_Y$ is a smooth del Pezzo surface and for any conic
$C\in\mathcal{C}^0_2(X)$
lying on
$S_Y$
one easily checks that the family of conics in $S_Y$ containing $C$
is a pencil containing 6 reducible conics. These reducible conics,
considered as pairs of lines, when $C$ runs through the open subset
$\mathcal{C}^{0*}_2(X):=v^{-1}(U)$ of $\mathcal{C}^0_2(X)$,
constitute the subset
$D_\kf^*$
of dimension 3, which is a dense open subset of $D_\kf$
since $D_\kf$ is irreducible of dimension 3 by Lemma \ref{Lines}.
Thus we have a map
$\psi:D_\kf^*\to\mathcal{C}^{0*}_2(X)$
which is 2:1 onto $D_1=\psi(D_\kf^*)$. Moreover, the map
$v$ factors as
$v:\mathcal{C}^{0*}_2(X)\tto{v_1}\kc\tto{v_2} U$
such that the maps $v_2$ and
$v_1|D_1$
are quasifinite and surjective,
$\deg(v_1|D_1)=6$ and $v_1^{-1}(c)\simeq\PP^1$ for $c\in\kc$.
Thus since
$D_\kf^*$
is irreducible, it follows that
$\mathcal{C}^{0*}_2(X)$ is irreducible of dimension 4.

Now take
$Y\in\check W$, so that
$Y$ is a tangent projective plane to $W$ at some point $x$:
\ $Y=PT_xW$.
Since $W$ is a general quartic in $\PP^3$ for a general $X$,
the set $T_Y$ of points $x$ such that
$Y=PT_xW$ is finite (generically it is a unique point).
Assume that
$\dim v^{-1}(Y)\ge2$.
Since
$\pi:S_Y\to Y$
is a finite morphism, the map
$\widetilde{v}:v^{-1}(Y)\to|\ko_Y(2)|:C\mapsto\pi(C)$
is quasifinite, hence
$\dim V(Y)\ge2$,
where
$V(Y):=\widetilde{v}(v^{-1}(Y))$.
Let
$V_1(Y)=\{C\in V(Y)|C\ni x$ for some $x\in T_Y\}$.
If $\dim V_1(Y)\ge2$,
then some point
$x\in T_Y$
belongs to at least a 1-dimensional family of reducible conics from
$V_1(Y),$
i.e. there is a 1-dimensional family of lines through $x$ in $Y$, each
of which is a double tangent to the quartic curve $W\cap Y$, which is
impossible when $W$ is general. Hence
$\dim V_1(Y)\le1$,
so that
$\dim v^{-1}(Y)'\ge2$,
where
$v^{-1}(Y)'=v^{-1}(Y)\smallsetminus\widetilde{v}^{-1}(V_1(Y))$.
Now for any conic
$C\in v^{-1}(Y)'$
the surface
$S_Y=\pi^{-1}(Y)$ is smooth along $C$,
so that, as above, the conic $C$ varies in a pencil of conics on $S_Y$,
contrary to the above inequality
$\dim_{C}v^{-1}(Y)'\ge2$.
Hence
$\dim v^{-1}(Y)\le1$ for $Y\in\check W,$
and so
$\dim v^{-1}(\check W)\le3.$ Now, by Riemann-Roch
$\chi(N_{C/X})=4$
for any conic
$C\in\mathcal{C}^0_2(X)$,
hence by deformation theory
$\dim_{C}\mathcal{C}^0_2(X)$
and the inequality
$\dim v^{-1}(\check W)\le3$
implies that
$v^{-1}(\check W)$
is not a component of
$\mathcal{C}^0_2(X)$.
Since by construction
$\mathcal{C}^0_2(X)=\mathcal{C}^{0*}_2(X)\cup v^{-1}(\check W)$,
the irreducibility of
$\mathcal{C}^0_2(X)$ follows.

(ii) Let $(l_1,l_2)\in D_\kf$ be generic. Then $l_1$ meets $l_2$
quasi-transversely (that is, with different tangents)
at some point $x\in X$ and, by Lemma \ref{Lines},
$
N_{l_i/X}\simeq 2\mathcal{O}_{{\PP}^1}.
$
Let $C=l_1\cup l_2$. We have the exact triples
$0\to N_{l_i/X}\to N_{C/X}|l_i\xrightarrow{\varepsilon_i}\mathbf{C}_x\to0,\ \ i=1,2,\ $
so that
$
N_{C/X}|l_i\simeq\mathcal{O}_{{\PP}^1}\oplus\mathcal{O}_{{\PP}^1}(1),\ \ i=1,2,\
$
and the exact triple
$0\to N_{C/X}\to N_{C/X}|l_1\oplus N_{C/X}|l_2\to N_{C/X}|_x\to0,\ \
N_{C/X}|_x\simeq\mathbf{C}_x^2,$
implies that
$
h^1(N_{C/X})=0.
$
Also,
the equality
$
N_{C/X}|l_1=\mathcal{O}_{{\PP}^1}\oplus\mathcal{O}_{{\PP}^1}(1),\
$
and the exact triple
$
0\to N_{C/X}|l_2(-x)\to N_{C/X}\xrightarrow{R}N_{C/X}|l_1\to 0,
$
in which $R$ is the usual restriction map, imply the surjectivity of the map
\mbox{$
H^0R:H^0(N_{C/X})\to H^0(N_{C/X}|l_1).
$}
Hence the composition
$
\tau_x:H^0(N_{C/X})\xrightarrow{H^0R}H^0(N_{C/X}|l_1)
\xrightarrow{H^0\varepsilon_1}H^0(\mathbf{C}_x)
$
is surjective. Hence
the equality
$
h^1(N_{C/X})=0
$
and \cite[Prop. 1.1]{HH} imply that
the conic $C=l_1\cup l_2$ is strongly smoothable in $X$,
that is there exists an analytic family
$\pi:\ks\overset{i}{\hookrightarrow}X\times\Delta\xrightarrow{pr_2}\Delta$
over the disk $\Delta$ such that $\ks$ is a smooth complex
surface,  $\pi$ is a flat morphism, $\pi^{-1}(0)=C$ and
$C_t=i(\pi^{-1}(t))\simeq\pi^{-1}(t)$
is a smooth conic in $X$ for any
$t\in\Delta\smallsetminus\{0\}$. Hence we can take a section
$s:\Delta\hookrightarrow\ks ,t\mapsto y_t,$ of $\pi$ such that
$y_0\ne x, y_0\in l_1$.
It follows now from
the above description of $N_{C/X}|l_1$
and the exact triple
$
0\to N_{C/X}|l_2(-x)\to N_{C/X}(-y_0)\to N_{C/X}|l_1(-y_0)\to0
$
that $h^1(N_{C/X})(-y_0)=0$, so that by semicontinuity
(possibly after shrinking
$\Delta$) we obtain $h^1(N_{C_t/X})(-y_t)=0$ for
$t\in\Delta\smallsetminus\{0\}$. Combining this with the observation that
$\det(N_{C_t/X})\simeq\ko_{{\PP}^1}(2)$ for
$\ \ t\in\Delta\smallsetminus\{0\}$,
we see that (ii) holds for a general conic.
\end{proof}
\end{sub}

\begin{sub}\label{Cubics} \rm {\bf Twisted cubics.}

\begin{lemma}\label{cubics}
Let $X$ be a general quartic double solid. Then the following assertions hold:

(i) $\mathcal{C}^0_3(X)$ is irreducible of dimension $6$.

(ii) For a general smooth cubic $C\in\mathcal{C}^0_3(X)$, the normal bundle
$N_{C/X}$ is isomorphic to
$\mathcal{O}_{{\PP}^1}(k)\oplus\mathcal{O}_{{\PP}^1}(4-k),\ 0\le k\le4$.
\end{lemma}

\begin{proof}
The proof is similar to that of Lemma \ref{conics}.
Let us start by proving (ii).
Using Lemmas \ref{Lines} and \ref{conics}, pick a general
line $l\in\kf$ and a conic
$C'\in\kc_2^0(X)$ meeting each other quasi-transversely
in a unique point $x$ and such that
$
N_{l/X}\simeq 2\mathcal{O}_{{\PP}^1},\
N_{C'/X}\simeq\mathcal{O}_{{\PP}^1}(k)\oplus\mathcal{O}_{{\PP}^1}(2-k),\ \ \ 0\le k\le1.
$
Let $C=l\cup C'$. Then we have the exact sequences
$
0\to N_{l/X}\to N_{C/X}|l\xrightarrow{\varepsilon}\mathbf{C}_x\to0,\
0\to N_{C'/X}\to N_{C/X}|C'\to\mathbf{C}_x\to0.
$
Hence
$
N_{C/X}|C'\simeq\mathcal{O}_{{\PP}^1}(a)\oplus\mathcal{O}_{{\PP}^1}(b),
\ 0\le a\le b,\ \ a+b=3,
$
and the exact triple
$
0\to N_{C/X}\to N_{C/X}|l\oplus N_{C/X}|C'\to N_{C/X}|_x\to0,\
N_{C/X}|_x\simeq\mathbf{C}_x^2,
$
implies that
$
h^1(N_{C/X})=0.
$
As in the proof of Lemma \ref{conics} (ii)
with  $(l,C')$ in place of  $(l_1,l_2)$,
we obtain
the wanted assertion
for a generic cubic from any irreducible component
of $\kc^0_3(X)$ which contains smoothings of a generic decomposable curve.
Let us show that there is only one such component.

By \cite{We}, Corollary (5.4), the curve $D_l\subset \kf$ of lines
in $X$ meeting a given generic line $l$ is irreducible and
is a 6-sheeted covering of $l$, because there are 6 lines passing
through the generic point of $X$. So, if we take a generic reducible
conic $C=l_1+l_2\subset X$, then the curve $D_C$ of lines
meeting $C$ is a 6-sheeted covering of $C$, which is irreducible
as a covering, that is there is no subcurve $D'\subset D_C$
which is a $d$-sheeted covering of $C$ for some $d<6$. The irreducibility
of a covering is stable under small deformations, so we
can conclude that $D_C$ stays irreducible when we deform $C$
into a smooth conic. The family of conics in $X$ being irreducible,
the incidence variety $\kd$ of pairs $(l,C')\in \kf\times\kc_2^0$
such that $l\cap C'\neq \empt$ is irreducible. By the above,
$l+C'$ is a smooth point of $\Hilb (X)$ for generic $(l,C')\in \kd$,
so there is only
one component $\kc^{0*}_3(X)$ of $\kc^{0}_3(X)$ containing the smoothings
of the curves $l+C'$ for generic $(l,C')\in \kd$.

Now we can prove the assertion (i) which implies, in particular,
that $\kc^{0*}_3(X)$ is the whole of  $\kc^{0}_3(X)$ and
hence (ii) is verified for a generic cubic $C$ from $\kc^{0}_3(X)$.

We have
$
\dim\kc^{0*}_3(X)=6,
$
and any irreducible component $\kc^{0'}_3(X)$ of $\kc^0_3(X)$ has dimension
at least $\chi (N_{C/X})=6$.

Let
$
\mathbf{X}\hookrightarrow Y\times{\PP}^{45}
$
be the total space of the linear system $\LLL^2$ on $Y$ and
$Y\xleftarrow{p}\mathbf{X}\xrightarrow{q}{\PP}^{45}$ the
natural projections.
Let
$
U=\{t\in{\PP}^{45}~|~X_t=p(q^{-1}(t))\ {\rm is\ smooth}\} .
$
Any $C\in \kc_3^0(X_t)$ for $t\in U$ is projected by
$\rho$ to a twisted cubic $C_0\subset{\PP}^{3}$.
Clearly,
$(\ko_{\PP^3}\oplus\ko_{\PP^3}(2))|_{C_0}\simeq
\ko_{\PP^1}\oplus\ko_{\PP^1}(6)$ for any
$C_0\in\kc_3^0(\PP^3)$.
Hence there exists a rank-8 vector bundle
$E$ over $\kc_3^0(\PP^3)$ with fiber
$E_{\left\{C_0\right\}}=H^0(\ko_{\PP^1}
\oplus\ko_{\PP^1}(6))$
over
$C_0\in\kc_3^0(\PP^3)$,
and
$ \kc_3^0(Y)$ is a dense open subset of
${\Proj}(E)$ such that the natural projection
$\kc_3^0(Y)\to\kc_3^0(\PP^3),\ C\mapsto\rho(C)$, coincides with the composition
$ \kc_3^0(Y)\hookrightarrow{\Proj}(E)\xrightarrow{\pr}\kc_3^0(\PP^3)$,
where $\pr$ is the structure morphism. $\kc_3^0(\PP^3)$ being
irreducible, $\kc_3^0(Y)$ is irreducible too and,
since
$\dim\kc_3^0(\PP^3)=12$,
we obtain
$\dim \kc_3^0(Y)=19.$

Consider the incidence variety
$
\Gamma=\{(t,C)\in U\times \kc_3^0(Y)~|~C\subset X_t\},\ \
U\xleftarrow{p_1}\Gamma\xrightarrow{p_2} \kc_3^0(Y).
$
We have
$
p_2^{-1}(C)=U\cap|\ki_{C/Y}\otimes\kl^{\otimes2}|
\overset{open}{\hookrightarrow}|\ki_{C/Y}\otimes\kl^{\otimes2}|.
$
Any $C_0\in\kc_3^0(\PP^3)$ has the following
resolution of its structure sheaf:
$0\to2\ko_{\PP^3}(-3)\to3\ko_{\PP^3}(-2)\to\ko_{\PP^3}\to\ko_{C_0}
\to0$.
Hence the twisted ideal sheaf
$\ki_{C/Y}\otimes\kl^{\otimes2}$ of any $C\in \kc_3^0(Y)$
has a resolution of the form
$
0\to2\kl(-3h)\to3\kl(-2h)\oplus2\kl^{\otimes2}(-3h)\to\kl\oplus3\kl^{\otimes2}(-2h)\to
\ki_{C/Y}\otimes\kl^{\otimes2}\to0.
$
An easy computation gives
$
|\ki_{C/Y}\otimes\kl^{\otimes2}|=\PP^{32}.
$
Hence $\Gamma$ is an irreducible variety of dimension 51.
As
$p_1:\Gamma\to U$ is dominant and $\dim U=45$, all the irreducible
components of
$
p_1^{-1}(t)= \kc_3^0(X_t),\  t\in U,
$
are generically of dimension 6. There is a distinguished one,
$ \kc_3^0(X_t)\cap\kc_3^{0*}(X_t)$, so by the irreducibility of
$\Gamma$ and $U$, \
$\kc_3^0(X_t)=\kc_3^{0*}(X_t)$.
\end{proof}
\end{sub}

\begin{sub}\label{Quintics} \rm {\bf Elliptic quintics.}
The following theorem is the main result of Section \ref{cocuqu}.
\begin{theorem}\label{quintics}
Let $X$ be general. Then $\kc^1_5(X)$
is irreducible and of dimension $10$.
\end{theorem}

\begin{proof} The proof is similar to that of Lemma
\ref{cubics}.

Pick a general cubic $C_1\in\kc_3^0$ and a smooth elliptic
conic $C_2\in\kc_2^1(X)$ meeting each other quasi-transversely
at a unique point $x$ and such that
$
N_{C_1/X}\simeq\mathcal{O}_{\PP^1}(k)\oplus\mathcal{O}_{\PP^1}(4-k),\ \
0\le k\le2,\ \
N_{C_2/X}\simeq 2\mathcal{O}_{C_2}(x+y),\ \ \{x,y\}=\pi^{-1}(\pi(x)).
$
Let $C=C_1\cup C_2$. Then we have the exact sequences
$
0\to N_{C_i/X}\to N_{C/X}|C_i\xrightarrow{\varepsilon_i}\mathbf{C}_x\to0,\ \  i=1,2.
$
Repeating the arguments from
the proof of Lemma \ref{cubics},
we see that there exists a
unique irreducible component
$\kc^{1*}_5(X)$
of $\Hilb (X)$ which contains the
generic reducible curve $C=C_1\cup C_2$ as above.
From the values of
$h^i(\NNN_{C/X})$, we see that
$
\dim\kc^{1*}_5(X)=10
$ and $
\dim\kc^{1'}_5(X)\ge10
$
for any other irreducible component $\kc^{1'}_5(X)$ of
$\kc^1_5(X)$.

For $t\in U$ any quintic $C\in \kc_5^1(X_t)$ is projected by
$\rho$ to an elliptic quintic $C_0\subset\PP^{3}$.
It is well known that  $\kc_5^1(\PP^3)$ is irreducible of dimension
20 \cite{Hu}. Arguing as in the proof of Lemma
\ref{cubics}, we obtain that
there exists a rank-11 vector bundle $E$ over $\kc_5^1(\PP^3)$
with fiber $E_{\left\{C_0\right\}}=H^0((\ko_{\PP^3}\oplus\ko_{\PP^3}(2))|_{C_0})$
over $C_0\in\kc_5^1(\PP^3)$, and
$ \kc_5^1(Y)$ is a dense open subset
of ${\Proj}(E)$ such that the natural projection
$ \kc_5^1(Y)\to\kc_5^1(\PP^3),\ C\mapsto\rho(C)$,
coincides with the composition
$ \kc_5^1(Y)\hookrightarrow{\Proj}(E)\xrightarrow{\pr}\kc_5^1(\PP^3)$,
where $\pr$ is the structure morphism.
Hence $ \kc_5^1(Y)$ is irreducible and
$
\dim \kc_5^1(Y)=30.
$

Consider the incidence variety
$
\Gamma=\{(t,C)\in U\times \kc_5^1(Y)~|~C\subset X_t\},\ \
U\xleftarrow{p_1}\Gamma\xrightarrow{p_2} \kc_5^1(Y).
$
We have
$ p_2^{-1}(C)=U\cap|\ki_{C/Y}\otimes\kl^{\otimes2}|
\overset{open}{\hookrightarrow}|\ki_{C/Y}\otimes\kl^{\otimes2}|.
$
Now take any quintic $C\in \kc_5^1(Y)$. Let
$\rho_C=\rho |_C:C\to C_0$ be the natural isomorphic projection of $C$
onto a quintic $C_0\in\kc_5^1(\PP^3)$. We have
$
p_{Y*}(\kl^{\otimes2}\otimes\ko_{C})\simeq\ko_{\PP^3}(4)|C_0.
$
Consider the standard exact triple
$
0\to\ki_{C/Y}\otimes\kl^{\otimes2}\to\kl^{\otimes2}\xrightarrow{\varepsilon}
\kl^{\otimes2}\otimes\ko_{C}\to0.
$
Applying $\rho_*$, we obtain the following one:
$
0\to p_{Y*}\ki_{C/Y}\otimes\kl^{\otimes2}\to
\ko_{\PP^3}\oplus\ko_{\PP^3}(2)\oplus\ko_{\PP^3}(4)\xrightarrow{\varepsilon_0}
\ko_{\PP^3}(4)|C_0\to0.
$
Here the map
$
\varepsilon_0|\ko_{\PP^3}(4):\ko_{\PP^3}(4)\to\ko_{\PP^3}(4)|C_0
$
coincides with the restriction map in the triple
$
0\to\ki_{C_0,\PP^3}(4)\to\ko_{\PP^3}(4)\xrightarrow{\res}\ko_{\PP^3}(4)|C_0\to0.
$
Hence,
by the above triples,
the equality
$
h^1(\ki_{C/Y}\otimes\kl^{\otimes2})=0
$
is a consequence of
$
h^1(\ki_{C_0,\PP^3}(4))=0.
$
To prove the latter, take two general cubic surfaces $S_1,S_2$
through $C_0$ such
that the residual curve of their intersection is a
smooth rational quartic $C_1$
meeting $C_0$ in a scheme $Z$ of length 10.
Then we have two exact triples
$
0\to\ki_{C_0\cup C_1,\PP^3}(4)\to\ki_{C_0,\PP^3}(4)\to\ko_{C_1}(4)(-Z)\to0,
$
where $\ko_{C_1}(4)(-Z)=\ko_{\PP^1}(6)$, and
$
0\to\ko_{\PP^3}(-2)\to2\ko_{\PP^3}(1)\to\ki_{C_0\cup C_1,\PP^3}(4)\to0.
$
Their long cohomology exact sequences imply the vanishing
of $h^1(\ki_{C_0,\PP^3}(4))=0$.
Hence $h^0(\ki_{C/Y}\otimes\kl^{\otimes2})=26$
and
$
|\ki_{C/Y}\otimes\kl^{\otimes2}|=\PP^{25}.
$
The remaining part of the proof is the same as in Lemma \ref{cubics}.
\end{proof}

In Section \ref{abjac}, we will use the following result
which can be easily proved by the techniques of \cite{HH}.
\begin{lemma}\label{Decomp}
Let $X$ be general so that $\CCC^1_5(X)$ is irreducible and
let $\overline{\CCC^1_5(X)}$ be the
closure of $\CCC^1_5(X)$ in $\Hilb(X)$. Then there exists a strongly smoothable reduced
nodal curve $C_0\in\overline{\CCC^1_5(X)}$ of the form
$
C_0=l_1\cup...\cup l_5, \ \ \ l_i\in \FFF,\ \ i=1,...,5,
$
such that its only singularities are the $5$ distinct points
$
a_1=l_1\cap l_2, \ \ a_2=l_2\cap l_3, \ \ a_3=l_3\cap l_4, \ \ a_4=l_4\cap l_1, \ \
a_5=l_1\cap l_5,
$
and
$
N_{l_i/X}=2\OOO_{l_i}, \ i=1,...,5,
$
$
N_{C_0/X}|l_1\simeq\OOO_{l_1}(1)\oplus\OOO_{l_1}(2),\ \
N_{C_0/X}|l_i\simeq2\OOO_{l_i}(1),\ 2\le i\le4, \ \
N_{C_0/X}|l_5\simeq\OOO_{l_5}\oplus\OOO_{l_5}(1).
$
\end{lemma}

\end{sub}


\section{The Abel-Jacobi map of the family of elliptic quintics}
\label{abjac}

In this section we will study the Abel-Jacobi maps of three families
of curves on a general double solid $X$:
genus-4 septics $\CCC^4_7(X)$, elliptic quintics $\CCC^1_5(X)$ and
a certain family $T$ of reducible curves of degree 17 and arithmetic genus
29 which will set a connection between the two other families.

\begin{sub}\rm {\bf Notation.}
We will use the following symbols, in addition to those introduced
in~\ref{nota-irred}:

\begin{itemize}
\item $\CCC^4_7(\PP^3)\subset \Hilb(\PP^3)$,
the base of the family of smooth irreducible septics of genus 4 in $\PP^3.$
\item $\CCC^4_7(X)$ (resp. $\CCC^4_7(Y)$),
the base of the family of
septics of genus 4 in $X$ (resp. $Y$) which are mapped by $\rho$
isomorphically onto a curve from $\CCC^4_7(\PP^3)$.
\item $\overline{\CCC^4_7(X)}$, the closure of $\CCC^4_7(X)$ in $\Hilb(X).$
\item $\overline{\CCC^4_7(X)}^*:=\{C\in\overline{\CCC^4_7(X)}~|~C$ is a reduced
curve with at worst
ordinary double points$\}$, a dense open subset of $\overline{\CCC^4_7(X)}$ containing
$\CCC^4_7(X)$.
\item $\CCC^2_5(\PP^3):=
\{C\in\Hilb(\PP^3)|\ C$
is a reduced quintic curve of arithmetic genus 2 and bidegree (3,2)
on a smooth quadric
surface $Q\subset\PP^3\}.$
\item $\pi^*\CCC^2_5(\PP^3):=
\{C\in\Hilb(X)|\ C=\pi^{-1}(\underline{C}),
\ \underline{C}\in\CCC^2_5(\PP^3)\}.$
\item $\CCC_{(3,3)}(\PP^3):=
\{C\in\Hilb(\PP^3)|\ C$
is a curve of bidegree (3,3)
on a smooth quadric surface $Q\subset\PP^3$
(that is, $C$ has arithmetic genus 4 and is a complete
intersection of $Q$ and a cubic surface)$\}.$
\item $\pi^*\CCC_{(3,3)}(\PP^3):=
\{C\in\Hilb(X)|\ C=\pi^{-1}(\underline{C}),
\ \underline{C}\in\CCC_{(3,3)}(\PP^3)\}.$
\item $\overline{\CCC^1_5(X)}$, the closure of $\CCC^1_5(X)$ in $\Hilb(X)$.
\item $\overline{\CCC^1_5(X)}^*=\{C\in\overline{\CCC^1_5(X)}~|~C$ is a reduced curve with at most
ordinary double points$\}$, an open subset of $\overline{\CCC^1_5(X)}$.
\item $\CCC_{5+2}:=\{C\in\Hilb(X)~|~C$ is a reduced curve of the form $C=C'\cup C''$,
where $C'\in\overline{\CCC^1_5(X)}^*,C''\in\CCC^1_2(X)$ and $C'$ meets $C''$
quasi-transversely in 3 points$\}$.
\item $\check\CCC_{5+2}:=\{C'\cup C''\in\CCC_{5+2}~|~C'\in\CCC^1_5(X)\}$
\item $\Phi_\BBB:\BBB\to J(X)$,
the Abel-Jacobi map of a given subvariety $\BBB$ of the Hilbert scheme of curves on $X$
(defined uniquely up to a choice of a reference point in $\BBB$).
\end{itemize}
\end{sub}

\begin{sub}\rm   {\bf Septics of genus 4.}
We are interested in the family of septics
$\CCC^4_7(X)$ because of the following result
of Welters \cite[6.18]{We}:

\begin{proposition}\label{AJ-septics}
For $X$ general, the Abel-Jacobi map
$\Phi_{\CCC^4_7(X)}:\CCC^4_7(X)\to J(X)$
is dominant.
\end{proposition}

We begin with the following lemma.
\begin{lemma}\label{C(5+2)}
Let $X$ be a general quartic double solid. Then:

(i) $\check\CCC_{5+2}$ is irreducible of dimension $21$ and the natural projection
$p'':\check\CCC_{5+2}\to \CCC^1_5(X),\ C'\cup C''\mapsto C'$, is surjective.

(ii) Let $\overline{\check\CCC_{5+2}}$ be the closure of $\check\CCC_{5+2}$ in
$\CCC_{5+2}$. There exists a strongly smoothable curve
$C'_0\cup C''_0\in\overline{\check\CCC_{5+2}}\smallsetminus\check\CCC_{5+2}$
such that
$
C'_0=l_1\cup...\cup l_5\in\overline{\CCC^1_5(X)}\smallsetminus\CCC^1_5(X)
$
is the curve
from Lemma \ref{Decomp}
and $C''_0\in\CCC^1_2(X)$ is a smooth elliptic conic such that
$
Z_b:=C'_0\cap C''_0=\{b_1,b_2,b_3\},\ \ \ b_1\in l_2,\ \ \ b_2\in l_4,\ \ \
b_3\in l_5,\ \ \ Z_b\cap Z_a=\empt,
$
where $Z_a={\rm Sing}C'_0$.

(iii) The set
$
\CCC_{5+2}^*:=\{C\in\check\CCC_{5+2}~|~C~
\mbox{\em is~strongly~smoothable}\}\subset\CCC_{5+2}\cap\overline{\CCC^4_7(X)}^*
$
is an irreducible dense open subset of $\CCC_{5+2}$
and the projection
$p'':\CCC_{5+2}^*\to \CCC^1_5(X)$, \mbox{$C'\cup C''\mapsto C'$},
is dominant.

(iv) Let $\CCC_7^{4*}(X)$ be an irreducible component
of $\overline{\CCC_7^4(X)}$
containing $\CCC_{5+2}^*$. Then
$\dim\CCC_7^{4*}(X)=14$
and
$\CCC_7^{4*}(X)^0=\{C\in\CCC_7^{4*}(X)~|~
C~\mbox{\em is~smooth~and}~\pi|C:C\to\pi(C)~\mbox{\em is~an~isomorphism}\}$
is a nonempty open subset of $\CCC_7^{4*}(X)$.
\end{lemma}
\begin{proof}
(i) $B=\pi (C')$
is an elliptic quintic in $\PP^3$.
For any $x\in B$, there are at most two
trisecant lines of $B$ passing through $x$, and the family
$Z$ of trisecant lines of $B$ is an irreducible curve.  To prove this,
consider the projection $\pi_x$ with center
$x$ from $B$  to $\PP^2$. The image of $B$ is a plane quartic
of geometric genus 1, so we have the following two
possibilities for its singular locus: 1) two singular points, each one
of which is either a node, or a cusp, 2) one tacnode.
The inverse image in $B$ of a singular point of $\pi_x(B)$ is
a pair of points, say $z_1,z_2$, and the line
$\overline{z_1z_2}$ is a trisecant of $B$, which meets
$B$ in three distinct points $x,z_1,z_2$. Hence
the family $Z$ of trisecant lines of $B$ is a curve.
To prove its irreducibility,
consider the incidence curve
$D=\{ (z,l)\in B\times Z\mid z\in l\}$ in $B\times Z$
with its natural
projections $B\xleftarrow{p}D\xrightarrow{q}Z$. The above argument
shows that $p$ is of degree 2. To see that both $D$ and $Z$ are
irreducible, it suffices to
verify that $p$ has at least one simple ramification point.
The latter is easily done in the manner of \cite{SR}, Chapter XII,
by interpreting such ramification points as fixed points of an appropriate
correspondence $G\subset B\times B$ and applying to $G$ the
Chasles-Cayley-Brill formula (\cite{SR}, Chapter XII, Corollary 3.52);
we leave the details to the reader. Now
${p''}^{-1}(C')$ is
identified with a nonempty open subset of $Z$, so it is 1-dimensional and
irreducible, which implies (i).

(ii) Remark that, in Lemma \ref{Decomp}, one can choose the curve
$C'_0=l_1\cup...\cup l_5$ in such a way that the lines
$m_i=\pi(l_i),\ i=2,4,5$ in $\PP^3$
are disjoint. Then it is obvious that a small
(in the classical topology over $\CC$)
smooth deformation $C'_t$ of $C'_0$ can be lifted to a small deformation
$C'_t\cup C''_t\in\check\CCC_{5+2}$ of $C'_0\cup C''_0$.
By the strong smoothability
of $C'_0$
(Lemma \ref{Decomp}),
$C:=C'_0\cup C''_0\in\overline{\check\CCC_{5+2}}$.
To see that $C$ is strongly smoothable, we consider the natural exact triples
: (a)
$
0\to N_a\to N_{C/X}\to\OOO_{Z_a}\to0,\ \ \
0\to N'_a\to N_a\to N_{C_1/X}\to0,
$
and (b)
$
0\to N_b\to N_{C_0/X}\to\OOO_{Z_b}\to0,\ \ \
0\to N'_b\to N_b\to N_{C''_0/X}\to0,
$
where
$C_1:=l_1\cup l_2\cup l_4\cup C'',\ \
N'_a:=(N_{C/X}|l_3\cup l_5)(-Z_a),\ \
N'_b:=(N_{C/X}|C'_0)(-Z_b)$.
Then the strong smoothability of $C$ will follow from the equalities
$
h^1(N_a)=0,
$
$
h^1(N_b)=0.
$
To prove the first of them, one deduces from Lemma \ref{Decomp}
that
$N'_a=2\OOO_{l_3}(-1)\oplus2\OOO_{l_5}(-1)$. Then the exact triples
$
0\to2\OOO_{l_2}(-1)\oplus2\OOO_{l_4}(-1)\to N_{C_1/X}\to N_{C_1/X}|C''_0\oplus
2\OOO_{l_1}(1)\to0,\ \ \
0\to N_{C''_0/X}\to N_{C_1/X}|C''_0\to\CC_{b_1}\oplus\CC_{b_2}\to0
$
together with the second triple
in (b)
give
$
h^1(N_a)=0.
$

To prove
that
$
h^1(N_b)=0,
$
we include $N'_b$ into the series of exact triples
following from the construction of $C$:
$
0\to N_1\to N'_b\to\OOO_{l_5}\oplus\OOO_{l_5}(-1)\to0,\ \ \
0\to2\OOO_{l_3}(-1)\to N_1\to N_2\to0,\ \ \
0\to N_3\to N_2\to\OOO_{l_2}\oplus\OOO_{l_2}(1)\to0,\ \ \
0\to\OOO_{l_4}\oplus\OOO_{l_4}(-1)\to N_3\to\OOO_{l_4}\oplus\OOO_{l_4}(-1)\to0.
$
As $C''_0=\pi^{-1}(m),$ where $m$ denotes a line in $\PP^3$, one immediately
verifies that $h^1(N_{C''_0/X})=0$. This together with the second triple
in (b)
gives
$
h^1(N_b)=0.
$

(iii) follows from (i) and (ii) since the strong smoothability is an open condition
in $\overline{\check\CCC_{5+2}}$.

(iv) Let
$C=C'_0\cup C''_0\in\overline{\check\CCC_{5+2}}\smallsetminus\check\CCC_{5+2}$
from (ii) above. From the above exact sequences, one immediately obtains
$h^1(N_{C/X})=0,\ h^0(N_{C/X})=14,$ hence $\dim\CCC_7^{4*}(X)=14$.
Next, consider $C$ as the union $C=C_3\cup C_4$, where
$C_3:=l_5\cup C''_0$ and $C_4:=l_1\cup l_2\cup l_3\cup l_4$. One easily
verifies that $C_3$ is strongly smoothable
and
$
\dim_{\{C_3\}}\Hilb(X)=6.
$
Thus there exists a unique 6-dimensional
irreducible component $\CCC_3^{1*}(X)$ of
$\Hilb(X)$ containing $C_3$.
Denote $R=\{\tilde C\in\CCC_3^{1*}(X)~|~\tilde C$ is
reducible of the form $\tilde C=l\cup\pi^{-1}(m)$,
where $l$ is a line in $X$ and $m$
is a line in $\PP^3$\}. By construction, $C_3\in R$.
It follows easily that:
a) for any
$\tilde C\in\CCC_3^{1*}(X)\smallsetminus R$,\
$\pi|\tilde C:\tilde C\to\pi(\tilde C)$
is an isomorphism of
$\tilde C$ onto a plane cubic $\pi(\tilde C)$;
b) $\dim R=5$, and $\dim_{\{C_3\}}R(l_1,l_2,l_4)=2$, where
$R(l_1,l_2,l_4)=\{\tilde C\in R~|~\tilde C$ intersects
each one of the lines $l_1,l_2,l_4$
in one point$\}$.  Thus,
since
$
\dim_{\{C_3\}}\Hilb(X)=6,
$
there exists
$\tilde C\in\CCC_3^{1*}(X)\smallsetminus R$ intersecting each of the lines $l_1,l_2,l_4$
in one point, so that $\pi$ maps $\tilde C\cup C_4$
isomorphically onto its image. By construction,
$\tilde C\cup C_4\in\CCC_7^{4*}(X)$. The property of a curve in $X$
to be mapped isomorphically onto its image in $\PP^3$ is open and a general
curve in $\CCC_7^{4*}(X)$ is smooth as a deformation of the strongly smoothable
curve $C_3\cup C_4$, so $\CCC_7^{4*}(X)^0$ is open in
$\CCC_7^{4*}(X)$, and this ends the proof.
\end{proof}

\begin{lemma}\label{sept-irr}
For $X$ general, $\CCC^4_7(X)$ is irreducible of dimension $14$.
\end{lemma}

\begin{proof}
This is proved in the same way as Lemma
\ref{cubics} and
Theorem
\ref{quintics}. First, remark that $\CCC^4_7(\PP^3)$
is an irreducible variety of dimension 28 (see, e.g., \cite[Ch.4]{H}), and for any
$C_0\in\CCC^4_7(\PP^3)$ we have
$(\OOO_{\PP^3}\oplus\OOO_{\PP^3}(2))|_{C_0}\simeq
\OOO_{C_0}\oplus\OOO_{C_0}(D),\ \ \deg D=14,$ so that
$h^0((\OOO_{\PP^3}\oplus\OOO_{\PP^3}(2))|_{C_0})=12$.
This
implies that there exists a rank-12 vector bundle $E$ over $\CCC_7^4(\PP^3)$
with fiber $E_{\left\{C_0\right\}}=H^0((\OOO_{\PP^3}\oplus\OOO_{\PP^3}(2))|_{C_0})$
over $C_0\in\CCC_7^4(\PP^3)$.
Thus $ \CCC_7^4(Y)$ is a dense open subset
of ${\Proj}(E)$ such that the natural projection
$ \CCC_7^4(Y)\to\CCC_7^4(\PP^3),\ C\mapsto\rho(C)$, coincides with the composition
$ \CCC_7^4(Y)\hookrightarrow{\Proj}(E)\xrightarrow{\pr}\CCC_7^4(\PP^3)$,
where $\pr$ is the structure morphism.
Hence $\CCC_7^4(Y)$ is irreducible of dimension 39.

Now take any $C\in \CCC_7^4(Y)$
and compute $\dim |\ki_{C/Y}\otimes\kl^{\otimes2}|$.
As in the proof of Theorem \ref{quintics}, it suffices to show that
$h^1(\ki_{C_0/\PP^3}(4))=0$.
To prove this, take two general cubic surfaces $S_1,S_2$
through $C_0$ (such cubics exist since, by Riemann-Roch, $h^0(\OOO_{C_0}(3))=18$.)
There are two possibilities for $S_1$ and $S_2$:
1)~both $S_1$ and $S_2$ are irreducible, and
2) both $S_1$ and $S_2$ are reducible. We will consider
the case 1), the case 2) being treated similarly.

In this case $S_1\cap S_2=C_0\cup C'$, where $C'$ is a curve of degree 2.
We have the following two subcases:
1.a) $C'$ a plane conic, or 1.b) $C'$ is a pair of two skew lines,
possibly coincident, that is a double line on a smooth quadric.
Consider the scheme $Z=C_0\cap C'$. One easily verifies that
$\chi(\OOO_Z)=7$ in case 1.a) and $\chi(\OOO_Z)=8$ in case 1.b),
and that $\chi(\OOO_{Z\cap l})=4$ for any line
$l\subset C'$. Thus, in both cases one obtains
$h^1(\III_{Z,C'}(4))=0$.
Then the exact triples
$
0\to\III_{C_0\cup C'/\PP^3}(4)\to\III_{C_0/\PP^3}(4)\to\III_{Z/C'}(4)\to0
$
and
$
0\to\OOO_{\PP^3}(-2)\to2\OOO_{\PP^3}(1)\to\III_{C_0\cup C'/\PP^3}(4)\to0
$
imply that $h^1(\ki_{C_0/\PP^3}(4))=0$.

Hence $h^0(\III_{C/Y}\otimes\LLL^{\otimes2})=
\chi(\III_{C/Y}\otimes\LLL^{\otimes2})=21$
and $|\III_{C/Y}\otimes\LLL^{\otimes2}|=\PP^{20}$ for $C\in \CCC_7^4(Y)$.
Thus, considering the incidence variety
$
\Gamma=\{(t,C)\in U\times \CCC_7^4(Y)~|~C\subset X_t\}$ with its natural
projections $
 U\xleftarrow{p_1}\Gamma\xrightarrow{p_2} \CCC_7^4(Y),
$
we see that $p_2^{-1}(C)$,
is open in $|\III_{C/Y}\otimes\LLL^{\otimes2}|$,
so that
$\Gamma$
is irreducible of dimension 59.
The remaining part of the proof goes along the same lines as the end
of the proofs of
Lemma
\ref{cubics} and
Theorem
\ref{quintics}.
Namely, by Lemma \ref{C(5+2)} (iv),
for a general double solid $X\in U$, there is a distinguished component
$\CCC^{4*}_7(X)$ of $p_1^{-1}(X)=\CCC^4_7(X)$
which plays the role of
$\CCC^{0*}_3(X)$, $\CCC^{1*}_5(X)$ in the proofs of the above mentioned lemma
and theorem.
(Remark that, since $\dim\Gamma=59=14+45=\dim\CCC^{4*}_7(X)+\dim U$,
$p_1:\Gamma\to U$ is dominant.)
\end{proof}

\end{sub}

\begin{sub}\rm  {\bf Curves of degree 17 and arithmetic genus 29.}
In this subsection we construct a family
$\tilde\Gamma_0$
of reducible curves $C=C_1\cup C_2$ of degree 17 and arithmetic genus
29 on $X$, where $C_1$ is a general septic from $\CCC^4_7(X)$ and
$C_2\in\pi^*\CCC^2_5(\PP^3)$ intersects $C_1$ in 13 points,
and a subfamily
$\Sigma_0\subset\tilde\Gamma_0$ (see (\ref{subfamily}))
of curves of the form $C=C'\cup C''\cup C_2$, where
$C'$ is a general quintic from
$\CCC^1_5(X)$
and $C''=\pi^{-1}(m)$, where $m$ is a line in $\PP^3$ trisecant to
$\pi(C')$. These two families will be used later, see
(\ref{Sigma}) and Proposition \ref{AJ(T)}.

Let
$$
M=\left\{(Q,l,Z,Z_l)~\left|\
\begin{minipage}{9 truecm}$Q$ is a smooth quadric in
$\PP^3$, $l$ a line on Q, $Z=\{y_4,...,y_{10}\}$  a
7-uple of distinct points on $Q\smallsetminus l$,
$Z_l=\{y_1,y_2,y_3\}$  a triple of distinct points on $l$
\end{minipage}\right.\right\}
$$
and
$
M^*=\{(Q,l,Z,Z_l)\in M~|~\dim|\JJJ_{Z\cup Z_l/Q}(3,2)|=1\},
$
where we assume that $l$ is of bidegree (0,1) on $Q$.
Clearly, $M$ is irreducible and $M^*$ is nonempty, hence dense in $M$.
Next, let
$$
\Pi=\{(C,(Q,l,Z,Z_l))\in\CCC^1_5(X)\times M~|~\underline{C}\cap Q=Z\cup Z_l,~
{\rm where}~\underline{C}=\pi(C)\}.
$$
It has two natural projections:
$\CCC^1_5(X)\xleftarrow{p_\Pi}\Pi\xrightarrow{q_\Pi}M$.
It is easy to see that $\Pi$ is nonempty. Indeed, for any
$C\in\CCC^1_5(X)$,
consider the smooth quintic
$\underline{C}=\pi(C)$,
take any trisecant line $l$ of $\underline{C}$ and any
smooth quadric $Q$ through $l$ and let
$Z=C\cap(Q\smallsetminus l)$, $Z_l=\uC\cap l$;
then, by construction, $(C,(Q,l,Z,Z_l))\in\Pi$.
We can factor $p_\Pi$ as
\begin{equation}\label{p_Pi}
\begin{array}{c}
p_\Pi:\Pi\xrightarrow{p'}\CCC_{5+2}\xrightarrow{p''}\CCC^1_5(X),
\\
p':(C',(Q,l,Z,Z_l))\mapsto C'\cup \pi^{-1}(l), \ \
p'':C'\cup C''\mapsto C',\end{array}
\end{equation}
and our construction shows that $p'$ is surjective
with an irreducible and rational fiber of dimension 7.
Besides, for
$C'\in\CCC^1_5(X),$
\ ${p''}^{-1}(C')$
is isomorphic to a dense open subset of the 1-dimensional family
of trisecant lines of the smooth quintic
$\pi(C')$ in $\PP^3$.
This description of the fibres of $p'$ and $p''$
together with (\ref{p_Pi}) implies, in particular, that
$\Pi$ is irreducible for
general $X$.

Now show that
\begin{equation}\label{Pi*}
\Pi^*:= {p'}^{-1}(\CCC_{5+2}^*)\cap q_\Pi^{-1}(M^*)\ne\empt.
\end{equation}
Choose any $C\in \CCC^1_5(X)$ and a trisecant line $l$ of
$\underline{C}$ such that
$C\cup\pi^{-1}(l)\in\CCC_{5+2}^*$, where
$\underline{C}=\pi(C)$. Then
$Z_l=\underline{C}\cap l=\{x^0_1,x^0_2,x^0_3\}$ is a
triple of distinct points. Through each $x^0_i$ there is, in addition to $l$,
one more trisecant $l^i_0$ of $\underline{C},$ and we denote
$\{x^0_{2i+2},x^0_{2i+3}\}=l^0_i\cap\underline{C}\setminus\{ x^0_i\} ,
\ \ i=1,2,3.$
Since $\underline{C}$ has degree 5, it follows immediately that
$l^0_1,l^0_2,l^0_3$ are disjoint, hence
there exists a unique smooth quadric $Q_0$ containing
$l,l^0_1,l^0_2$ and $l^0_3$.
Let $x^0_{10}$
be the residual point of intersection of $\underline{C}$ and $Q^0$:
$\underline{C}\cap Q_0=\{x^0_1,...,x^0_{10}\}$.
One can easily verify that for sufficiently general $l$,
the 10 points $x^0_1,...,x^0_{10}$ are distinct.
Denote $Z^0=\{x^0_4,...,x^0_{10}\}$ and let
$l^0_4$
be the line on $Q^0$ through $x^0_{10}$ in the same ruling as $l$,
say $|\OOO_{Q_0}(0,1)|$. Then
$
|\JJJ_{Z^0\cup Z_l/Q_0}(3,2)|=\bigcup_{i=1}^4l^0_i+|\OOO_{Q_0}(0,1)|\simeq\PP^1.
$
Now consider the triple $x_0=(x^0_1,x^0_2,x^0_3)$ as a point of $\Sym^3(l)\simeq\PP^3$
and take a sufficiently small neighborhood $U$ of
$x_0$ in the classical topology of $\Sym^3(l)$
so that for any $x=(x_1,x_2,x_3)\in U$ the following properties are
satisfied: (i)
for each $i=1,2,3$, there is a unique trisecant line $l_i$
of $\underline{C}$ passing through $x_i$ and close to $l^0_i$
in the sense that
$\{x_{2i+2},x_{2i+3}\}=l_i\cap\underline{C}\setminus\{ x^0_i\}$
is close to the pair
$\{x^0_{2i+2},x^0_{2i+3}\}$ in $\Sym^2(C)$;
(ii) there exists a unique quadric $Q$
containing $l,l_1,l_2,l_3$ and close to $Q_0$ in the space of quadrics;
for this quadric, the residual point
$x_{10}=\underline{C}\cap Q\smallsetminus\{x_1,...,x_9\}$
is close to $x^0_{10}$.

Denote $Z_l=\{x_1,x_2,x_3\},\ \ Z=\{x_4,...,x_{10}\}$ and let
$l_4$ be the line on $Q$ through $x_{10}$
belonging to the same ruling as $l$.
These data define a point $m=(Q,l,Z,Z_l)\in M$.
Now we will determine the exact sense of the condition
that $U$ is sufficiently small:
\begin{equation}\label{1-dim}
|\JJJ_{Z\cup Z_l/Q}(3,2)|\simeq\PP^1.
\end{equation}
Remark that, by construction,
$\bigcup_{i=1}^4l_i+|\OOO_Q(0,1)|\subset|\JJJ_{Z\cup Z_l/Q}(3,2)|$,
so that
$
|\JJJ_{Z\cup Z_l/Q}(3,2)|=\bigcup_{i=1}^4l_i+|\OOO_Q(0,1)|.
$
The condition (\ref{1-dim}) is open and nonempty. So,
for any $x_0\in \Sym^3(l)$, we obtain a point
$m=m(x_0)=(Q,l,Z,Z_l)\in M^*$.
Then
$(C,m)\in\Pi^*$, and (\ref{Pi*}) follows.
As $\Pi$ is irreducible and $M^*$ is open in $M$,
this implies that $\Pi^*$ is dense and open in $\Pi$
and, consequently,
\begin{equation}\label{dom0}
p':\Pi^*\to\CCC_{5+2}^*\ \mbox{is dominant.}
\end{equation}

Hence (\ref{p_Pi}) and Lemma \ref{C(5+2)} (i) show that
\begin{equation}\label{dom1}
p_\Pi:\Pi^*\to\CCC^1_5(X),\ (C',(Q,l,Z,Z_l))\mapsto C',\ \mbox{is dominant.}
\end{equation}
Remark that, for any
$C_2\in \pi^*|\JJJ_{Z\cup Z_l/Q}(3,2)|$, we have :
\begin{equation}\label{3,3}
C''\cup C_2\in\pi^*\CCC_{(3,3)}(\PP^3).
\end{equation}

Define, for any given $z=(C',(Q,l,Z,Z_l))\in\Pi^*$,
\begin{equation}\label{C_1}\begin{array}{c}
C_1=p'(z)=C'\cup C'',\ \ \ C''=\pi^{-1}(l),
\\
f_z=\{C_2\in\pi^*\CCC^2_5(\PP^3)~|~C_2=\pi^{-1}(\underline{C}_2),\ \underline{C}_2\in
|\JJJ_{Z\cup Z_l/Q}(3,2)|\},\ f_z\simeq\PP^1,
\\
f_z^*=\{C_2\in f_z~|~C_2~\mbox{is~a~reduced~curve}\} ,\end{array}
\end{equation}
Clearly, $f_z^*$ is a nonempty open subset of $f_z$.
By construction, for any $C_2\in f_z^*$, the curves $C_1$ and $C_2$ are reduced and
meet each other quasi-transversely in 13 points:
\begin{equation}\label{intrsn1}
C_1\cap C_2=\{y_1,y_2,y_3,y'_1,y'_2,y'_3,y_4,...,y_{10}\},\ \ \
\tilde Z_l=\pi^{-1}(Z_l)=
\end{equation}
$$
=\{y_1,y_2,y_3,y'_1,y'_2,y'_3\}=C''\cap C_2,\
\tilde Z=(\pi|C)^{-1}(Z)=\{y_4,...,y_{10}\}=C'\cap C_2.
$$
Hence
\begin{equation}\label{C0}
C_0=C'\cup C''\cup C_2
\end{equation}
is a curve on $X$ of degree 17 and arithmetic genus 29.
Remark also that $\pi^*\CCC^2_5(\PP^3)\simeq\CCC^2_5(\PP^3)$
is a rational irreducible
variety of dimension~20.

Now consider the incidence variety
$$
\Gamma=\left\{(C_1,C_2)\in\overline{\CCC^4_7(X)}^*
\times\pi^*\CCC^2_5(\PP^3)~\left|~\begin{minipage}{7 truecm}$C_1$ and $C_2$
are reduced and meet each other quasi-transversely in 13 points
\end{minipage}\right.\right\}$$
with projections
$\overline{\CCC^4_7(X)}^*\xleftarrow{p_1}\Gamma\xrightarrow{p_2}
\pi^*\CCC^2_5(\PP^3)$
and define its locally closed subset
$
\Sigma^*=p_1^{-1}(p'(\Pi^*)).
$
From (\ref{p_Pi}), (\ref{dom0}) and (\ref{dom1}) it follows immediately that the natural
projection
\begin{equation}\label{p_Si*}
p_{\Sigma^*}:=p''\circ p_1:\Sigma^*\to \CCC^1_5(X),\
(C'\cup C'',C_2)\mapsto C',~\mbox{is dominant}.
\end{equation}
Moreover, for
$z\in\Pi^*,$   by (\ref{C_1})
and
the description of $p'$ in (\ref{p_Pi}), we
have
for $C_1=p'(z)$:
\begin{equation}\label{fib(p1)}
p_1^{-1}(C_1)~\mbox{is~a~rational~7-dimensional~variety}.
\end{equation}
Hence,
for general $X$, $\Sigma^*$ is
irreducible
of dimension 17, and there is an isomorphism
\begin{equation}\label{bij1}
\psi:\Sigma^*\xrightarrow{\sim}\Sigma' ,
(C'\cup C'',C_2)\mapsto C'\cup C''\cup C_2,
\end{equation}
where
$
\Sigma'=\{C_0\in\Hilb(X)~|~C_0~\mbox{is~defined~by~(\ref{intrsn1})-(\ref{C0})
for some $z\in\Pi^*$}\}.
$
In view of (\ref{bij1}), we rewrite (\ref{p_Si*}) as follows:
$$
p_{\Sigma'}:\Sigma'\to\CCC^1_5(X),\ C_0=(C'\cup C''\cup C_2)\mapsto C',~
\mbox{is dominant}.
$$
Also, by (\ref{3,3}), we have
$
\im(q_{\Sigma'})\subset\pi^*\CCC_{(3,3)}(\PP^3),
$
where $q_{\Sigma'}:C_0=(C'\cup C''\cup C_2)\mapsto C''\cup C_2$.

Next, looking at the local equations of $\Gamma$ in
$W=\overline{\CCC^4_7(X)}^*\times\pi^*\CCC^2_5(\PP^3)$ at a point
$
w=(C_1,C_2)\in\Sigma^*$,
we see that the transversal intersection of
$C_1$ and $C_2$ in each
one of the 13 points (\ref{intrsn1}) gives one local
equation for $\Gamma$ in $Y$; hence,
since
$
\dim\pi^*\CCC^2_5(\PP^3)=20
$, $\dim p_1^{-1}p_1(w)\geq
20-13=7$. Since by (\ref{fib(p1)}) this dimension is precisely 7,
it follows that
there exists an irreducible component of $\Gamma$, say $\Gamma_0$,
which contains $w$ and dominates $\overline{\CCC^4_7(X)}^*$
via $p_1$:
\begin{equation}\label{domin}
\overline{p_1(\Gamma_0)}=\overline{\CCC^4_7(X)}^*,
\end{equation}
where the closure is taken in $\overline{\CCC^4_7(X)}^*$
and we use the irreducibility of
$\overline{\CCC^4_7(X)}^*$ by Lemma \ref{sept-irr}. Hence,
taking into account the irreducibility
of $\Sigma^*$, (\ref{fib(p1)}) and
(\ref{bij1}), we see that
$
\Sigma_0=\psi(\Gamma_0\cap\Sigma^*)
$
is a dense open subset of $\Sigma$; also, by (\ref{p_Si*}),
the natural projection
\begin{equation}\label{p_Si0}
p_{\Sigma_0}:\Sigma_0\to\CCC^1_5(X),\ C_0=(C'\cup C''\cup C_2)\mapsto C',~
\mbox{is dominant}.
\end{equation}

Now consider in $\Hilb(X)$ the
subset
$
\tilde\Gamma_0=\{C\in\Hilb(X)
| C\ {\rm is\ a\ reducible\ curve\ of\ the}$
${\rm form\ }C=C_1\cup C_2,
\mbox{where $(C_1,C_2)\in\Gamma_0$}\}.
$
There is an isomorphism
$
\psi:\Gamma_0\xrightarrow{\sim}\tilde\Gamma_0,\ (C_1,C_2)\mapsto C_1\cup C_2,
$
similar to (\ref{bij1}), and by construction,
\begin{equation}\label{subfamily}
\Sigma_0\subset\tilde\Gamma_0.
\end{equation}
\end{sub}

\begin{sub}\rm    {\bf Criterion of Welters.}

Now we invoke the following important technical result of
G.Welters \cite[Prop. 6.17]{We}, giving a sufficient
condition for a family of curves to parametrize, via the
Abel-Jacobi map,  a translate of the theta divisor of $J(X)$.

\begin{proposition}\label{Welters}
Let $T$ be a smooth irreducible (not necessarily complete) variety
parametrizing a family $\{C_t|t\in T\}$ of reduced curves of degree
$d=k(k+1)-3$ on $X$ with the following properties:
(i) The Abel-Jacobi map $\Phi_T:T\to J(X)$ is dominant;
(ii) for all $t\in T$, the linear system of surfaces
$|\JJJ_{C_t/X}(k)|$
has dimension $k$
and the surfaces of this linear system do not contain all the lines of $X$;
(iii) the incidence divisor in $T\times \FFF$ is reduced.
Then for any line $l\in\FFF$,
the subvariety of $T$ defined by
$$
T_l=\{t\in T\mid C_t\cap l=\empt~\mbox{\em and}~C_t\cup
l~\mbox{\em lies~on~a~surface~from}~|\OOO_X(k)|\}
$$
is sent by the Abel-Jacobi map into a translate of the theta divisor
$\Theta \subset J(X)$.
\end{proposition}

Let us construct such a $T$ from the curves introduced
in the previous subsection.
Consider the Abel-Jacobi map
$\Phi_{\tilde\Gamma_0}:\tilde\Gamma_0\to J(X)$
and the induced map
\begin{equation}\label{Phi0}
\Phi_0:\Gamma_0
\isom{\psi}
\tilde\Gamma_0
\xrightarrow{\Phi_{\tilde\Gamma_0}}J(X),\ (C_1,C_2)
\mapsto\Phi_{\tilde\Gamma_0}(C_1\cup C_2),
\end{equation}
where
\begin{equation}\label{Phi0'}
\Phi_{\tilde\Gamma_0}(C_1\cup C_2)=
\Phi_{\overline{\CCC^4_7(X)}^*}(C_1)+\Phi_{\pi^*\CCC^2_5(\PP^3)}(C_2),\ \ \ \
C_1\in\overline{\CCC^4_7(X)}^*,\ \ C_2\in\pi^*\CCC^2_5(\PP^3).
\end{equation}
As $\pi^*\CCC^2_5(\PP^3)\simeq\CCC^2_5(\PP^3)$ is rational,
\begin{equation}\label{c}
\Phi_{\pi^*\CCC^2_5(\PP^3)}=\const,\ \
\im(\Phi_{\pi^*\CCC^2_5(\PP^3)})=\{c\} .
\end{equation}
Hence, taking the closure of the image of the Abel-Jacobi map in $J(X)$
and using (\ref{domin}),
(\ref{Phi0})-(\ref{c}) and Proposition~(\ref{AJ-septics}), we obtain:
\begin{multline}\label{domin2}
\overline{\im(\Phi_{\tilde\Gamma_0})}=\overline{\im(\Phi_0)}=
\overline{c+\im(\Phi_{\overline{\CCC^4_7(X)}^*}(\overline{p_1(\Gamma_0)}~))}
=\overline{c+\im(\Phi_{\overline{\CCC^4_7(X)}^*})}\\ =
c+\overline{\im(\Phi_{\CCC^4_7(X)})}=c+J(X)=J(X),
\end{multline}
that is the Abel-Jacobi map
$\Phi_{\tilde\Gamma_0}:\tilde\Gamma_0\to J(X)$
is dominant. Now, consider the Abel-Jacobi map
$\Phi_{\Sigma_0}:\Sigma_0\to J(X)$. As $\pi^*\CCC_{(3,3)}(\PP^3)$ is rational,
$\im(\Phi_{\pi^*\CCC_{(3,3)}(\PP^3)})=\{c'\}$ is a point,
and we obtain from (\ref{p_Si0}):
\begin{equation}\label{AJ(Si0)}
\overline{\im(\Phi_{\Sigma_0})}=\overline{\im(\Phi_{\CCC^1_5(X)})}+c'.
\end{equation}
Now let
$\sigma:T\xrightarrow{bir}\tilde\Gamma_0$
be any desingularization of $\tilde\Gamma_0$. Using (\ref{subfamily}), define
\begin{equation}\label{Sigma}
\Sigma\subset\sigma^{-1}(\Sigma_0),\ \mbox{an irreducible component of}\
\sigma^{-1}(\Sigma_0)\ \mbox{such that}\
\sigma:\Sigma\twoheadrightarrow\Sigma_0.
\end{equation}
Consider the Abel-Jacobi maps
$\Phi_{T}=\Phi_{\tilde\Gamma_0}\circ\sigma:T\to J(X)$
and
$\Phi_{\Sigma}=\Phi_{\Sigma_0}\circ\sigma:\Sigma\to J(X)$.
Then (\ref{AJ(Si0)}) gives
\begin{equation}\label{AJ(Si)}
\overline{\im(\Phi_\Sigma)}=\overline{\im(\Phi_{\CCC^1_5(X)})}+c'.
\end{equation}
We obtain the following proposition:

\begin{proposition}\label{AJ(T)}
For $X$ general, there exists a smooth irreducible
variety $T$ parametrizing some family of curves of degree $17$ and
arithmetic genus $29$ such that:
(i) the Abel-Jacobi map $\Phi_{T}:T\to J(X)$ is dominant;
(ii) $T$ contains the subset $\Sigma$, defined by (\ref{Sigma})
and satisfying (\ref{AJ(Si)}), such that the curves
$C_t$ for $t\in\Sigma$ split into three components
as in (\ref{intrsn1})-(\ref{C0});
(iii) the incidence divisor
$\DDD=\{(C,l)\in T\times\FFF~|~C\cap l\ne\empt\}$, endowed with its natural
scheme-theoretic structure, is reduced.
\end{proposition}
\begin{proof}
(i) and (ii) follow from (\ref{domin2}) and (\ref{AJ(Si)}).
Let us prove (iii). There are 6 distinct lines, counted
with multiplicity 1 in $\FFF$, which pass through
a general point of $X$. So, it is sufficient to show
that there exists a curve $C$
from $T$ which passes through a general point of $X$: in this case,
$\DDD_C$ is a reduced 6-sheeted covering of $C$, and hence $\DDD$
is reduced. But this is clear, since by the construction of $\Sigma$,
the components of the reduced curves $C_0\in\Sigma$
sweep all of $X$ as $C_0$ runs
through $\Sigma$.
\end{proof}

We proceed now to the proof of the remaining hypotheses
of the Welters criterion for our family $T$.

\begin{lemma}\label{h^0=4} Let $X$ be a general quartic double solid. Then
the following assertions are verified:
(i) For any
$C\in\Sigma$,
$h^0(\JJJ_{C/X}(4))=5;$
(ii) for a general (hence any) line
$l\in\FFF$
such that
$C\cap l=\empt$,
there exists a surface
$S\in|\OOO_X(4)|$
containing
$C\cup l$;
(iii) for general $C\in\Sigma$, the surfaces of the linear system
$|\JJJ_{C/X}(4)|$ do not contain all the lines of $X$.
\end{lemma}
\begin{proof}
(i) Let $C\in\Sigma$.
By the definition of $\Sigma$ we have
$C=C_1\cup C_2$,
where
$C_1\in\CCC^1_5(X)$,
$\pi:C_1\to \uC_1:=\pi(C_1)$
is an isomorphism,
$C_2=\pi^{-1}(\uC_2),$
where
$\uC_2\in\CCC_{(3,3)}(\PP^3)$, and
$Z:=C_1\cap C_2=\{y_1,...,y_{10}\}$
is a subset of the 13-uple of points (\ref{intrsn1}).
This means, in particular, that
$\uC_2$
lies in a uniquely defined quadric, say
$Q$. Let
$Z_0:=\pi(Z),\ \ \deg(Z_0)=\deg(Z)=10$.
Since by construction
$\uC_1\cap Q\supset Z_0$
(as schemes) and
$\deg(\uC_1\cap Q)=10=\deg(Z_0)$,
it follows that, scheme-theoretically,
$
Z_0=\uC_1\cap Q.
$
Using the fact that
$\OOO_Q(\uC_2)=\OOO_Q(3)$,
we get the following exact triple:
\begin{equation}\label{triple_1}
0\to\JJJ_{Q\cup \uC_1/\PP^3}(4)\to\JJJ_{\uC_1\cup \uC_2/\PP^3}(4)\to
\OOO_Q(1)\to0.
\end{equation}
Let $S\subset\PP^3$ be any smooth cubic surface
passing through $\uC_1$.
We have the exact triple
$
0\to\JJJ_{Q\cup S/\PP^3}\to\JJJ_{Q\cup \uC_2/\PP^3}\to
\JJJ_{(S\cap Q)\cup \uC_1/S}\to0.
$
Here
$\JJJ_{(S\cap Q)\cup \uC_1/S}\simeq\OOO_S(-2)(-\uC_1)$,
because $S$ is smooth, hence
twisting by
$\OOO_{\PP^3}(4)$,
we obtain:
$
0\to\OOO_{\PP^3}(-1)\to\JJJ_{Q\cup \uC_2/\PP^3}(4)\to
\OOO_S(2)(-\uC_1)\to0.
$
One easily verifies that the natural map
$H^0(\OOO_S(2))\to H^0(\OOO_{\uC_1}(2))$
is surjective, so that
$h^1(\OOO_S(2)(-\uC_1))=0$.
Besides,
$h^1(\OOO_{\PP^3}(-1)=0.$
Hence
the last triple
implies
$h^1(\JJJ_{Q\cup \uC_1/\PP^3}(4))=0.$
This together with (\ref{triple_1}) and the equality
$h^1(\OOO_Q(1))=0$
gives
$
h^1(\JJJ_{\uC_1\cup \uC_2/\PP^3}(4))=0.
$
Hence the exact triple
$0\to\JJJ_{\uC_1\cup \uC_2/\PP^3}(4)\to\JJJ_{\uC_2/\PP^3}(4)\to
\JJJ_{Z_0/\uC_1}(4)\to0$
yields the surjectivity
\begin{equation}\label{surj}
H^0(\JJJ_{\uC_2/\PP^3}(4))\twoheadrightarrow H^0(\JJJ_{Z_0/\uC_1}(4)).
\end{equation}
Now, since
$\pi:C_1\to \uC_1$
is an isomorphism, we have
$\pi_*\JJJ_{Z/C_1}(4)\simeq\JJJ_{Z_0/\uC_1}(4)$,
which provides the exact triple
$
0\to\JJJ_{Z_0/\uC_1}(4)\to\pi_*\OOO_C(4)\to\pi_*\OOO_{C_2}(4)\to0 .
$
Hence
$
h^0(\pi_*\OOO_C(4))=\chi(\OOO_C(4))=40.
$
Next, consider the natural restriction map
$g:\pi_*\OOO_X(4)\twoheadrightarrow
\pi_*\OOO_{C}(4)\twoheadrightarrow\pi_*\OOO_{C_2}(4)$.
Here
$\ker(g)=\pi_*\JJJ_{C_2/X}(4)$,
so passing to the map on global sections
$
f:H^0(\pi_*\OOO_X(4))\to H^0(\pi_*\OOO_C(4)),
$
we obtain:
$
{\rm coker}(f)={\rm coker}(H^0(\pi_*\JJJ_{C_2/X}(4))
\xrightarrow{\res}H^0(\JJJ_{Z_0/\uC_1}(4))).
$
By the projection formula for $\pi$, we have
$\pi_*\JJJ_{C_2/X}(4)=\pi_*\pi^*\JJJ_{\uC_2/\PP^3}(4)=
\JJJ_{\uC_2/\PP^3}(4)\oplus\JJJ_{\uC_2/\PP^3}(2)$,
and the composition
$H^0(\JJJ_{\uC_2/\PP^3}(4))\into H^0(\pi_*\JJJ_{C_2/X}(4))
\xrightarrow{\res}H^0(\JJJ_{Z_0/\uC_1}(4))$
coincides with the surjective map (\ref{surj}),
so
${\rm coker}(f)=0$.
Hence
the equalities
$h^0(\pi_*\OOO_C(4))=40$
and
$h^0(\OOO_X(4))=h^0(\OOO_{\PP^3}(4))\oplus h^0(\OOO_{\PP^3}(2))=45$
imply that
$
H^0(\JJJ_{C/X}(4))=H^0(\pi_*\JJJ_{C/X}(4))=\ker(f)=\CC^5.
$

(ii) Take a cubic surface $F$ containing
$\uC_1$ and write $\uC_2=F\cap Q,$ where $Q$ is a general quadric.
Then $C_1\cup C_2\subset
\pi^{-1}(F)$.
Take
a general plane $\PP^2$ in $\PP^3$ and
a line $l$ in the del Pezzo surface
$\pi^{-1}(\PP^2)$ which does not meet
$C_1\cup C_2$. Clearly such
$l$
exists for a general choice of $\PP^2$. Then
$S=\pi^{-1}(F\cup\PP^2)\in |\OOO_X(4)|$ is the desired surface.

(iii) This is a standard dimension count.
One embeds $C$  into
$Y$
and considers the linear system
$L=|\JJJ_{C/Y}\otimes\OOO_{Y/{\PP}^3}(2)|$.
From (i), we deduce that $\dim L=5$. Take now
a line $l\subset Y$, mapped by $\rho$ isomorphically
onto a line in $\PP^3$. If it is sufficiently general,
then there exists a unique
smooth $X\in L$ containing $l$; this $X$ is the wanted double solid, and
the lemma is proved.
\end{proof}

We are ready now to prove the main result of this section:

\begin{theorem}\label{Theta}
$\overline{{\rm Im}(\Phi_{\CCC^1_5(X)})}=\Theta+{\rm const}.$
\end{theorem}

\begin{proof}
The family
$T$ constructed in
Proposition \ref{AJ(T)} satisfies the hypotheses of
Proposition \ref{Welters} with $k=4$.
According to Lemma \ref{h^0=4}, for general  $l\in\FFF$,
we have $T_l\supset\Sigma$, where we replace $\Sigma$ by a
nonempty open subset, if necessary.
Hence Proposition \ref{Welters} and (\ref{AJ(Si)}) imply:
\begin{equation}\label{Qntcs-inThta}
\overline{\im(\Phi_{\CCC^1_5(X)})}\subset\Theta+{\rm const} .
\end{equation}
To prove the opposite inclusion,
let us recall a well known description of the image of the differential of
the Abel-Jacobi map (see, e.g., \cite[(2.11.0) and (6.19)]{We}),
in applying it specifically to the map
$\Phi:=\Phi_{\CCC^1_5(X)}:\CCC^1_5(X)\to J(X)$.
Take any unobstructed quintic
$C\in \CCC^1_5(X)$, so that
$T_C\CCC^1_5(X)=H^0(N_{C/X})=\CC^{10}$. Assume, moreover,
that $C\not\subset W$ and  denote
$C_0=\pi(C),\ \ Z_0=\pi(C\cap W)$, where $W$ is the branch quartic of $X$.
Consider the composition
$
\phi:\ \CC^{10}=T_C \CCC^1_5(X)\xrightarrow{d\Phi}T_{\Phi(C)}J(X)
\xrightarrow{\sim}T_0J(X)\simeq(H^0(\OOO_{\PP^3}(2)))^{\vee}.
$
Then the kernel of the dual map $\phi^{\vee}$ is described as follows:
$
\Ker(\phi^{\vee})= \Ker\{\CC^{10}=H^0(\OOO_{\PP^3}(2))\xrightarrow{H^0(r_{C_0})}
H^0(\OOO_{C_0}(2))\xrightarrow{H^0(r_{Z_0})}H^0(\OOO_{Z_0}(2))\},
$
where $r_{C_0}:\OOO_{\PP^3}(2)\to\OOO_{C_0}(2)$,
$r_{Z_0}:\OOO_{C_0}(2)\to\OOO_{Z_0}(2)$
are the restriction maps, and $\im\phi$ is the orthogonal complement of
$\Ker(\phi^{\vee})$.
By definition, $C_0$ is an elliptic quintic
in $\PP^3$, hence it does not lie on a quadric, so
$H^0(r_{C_0})$ is an isomorphism, and $\Ker(H^0(r_{Z_0}))\simeq H^0(\OOO_{C_0})
=\CC$,
as follows from the exact triple $0\to\OOO_{C_0}\xrightarrow{\cdot Z_0}\OOO_{C_0}(2)
\xrightarrow{r_{Z_0}}\OOO_{Z_0}(2)\to0$. Hence
$\dim\Ker(\phi^{\vee})=1$ and $\dim\im\Phi=9$.
This together with (\ref{Qntcs-inThta}) and the irreducibility of $\Theta$ (see
\cite{We}) implies the wanted assertion.
\end{proof}

\end{sub}

\section{Pentagons on a degenerate double solid}

\begin{sub}\rm {\bf Degeneration into a pair of ${\mathbf P}^
{\boldsymbol 3}$'s.} \label{deg-pairs}

In this section, we will use the technique of degeneration
of the double solid into a reducible variety, developed by
Clemens in \cite{C-1}, \cite{C-2}. Let $\ppi :\X\lra\Delta$ be
a family of double covers $\pi_t:X_t\lra \PP^3$
branched in the quartics $W_t=\{ tF+G^2=0\}\subset \PP^3$
($t\in\Delta$), where $\Delta$ is an open disc in $\CC$, the equation
$F=0$ defines a smooth quartic $W$ and
$G=0$ a smooth quadric $Q$ such that $W\cap Q=B$ is a smooth
octic curve. For $t\in\Delta^*=\Delta\setminus\{ 0\}$, the fiber $X_t$
is a smooth double solid, and for $t=0$ a union of two copies
of $\PP^3$, say $\PP ',\PP ''$, meeting each other transversely along $Q$.

One can associate to $\ppi$ the compactified family $\F\lra\Delta$
of Fano surfaces, whose fiber ${\FFF}_t$ over $t\neq 0$ parametrizes lines in $X_t$,
and the Neron model $\J\lra\Delta$ of the family of the
intermediate Jacobians $J_t$ of the threefolds $X_t$.
The compactifications of both families are described by Clemens.
Namely, ${\FFF}_0$ is the union of two components ${\FFF}'\cup {\FFF}''$, where
${\FFF}'$ (resp. ${\FFF}''$) parametrizes the bisecant lines to $C$ in
$\PP '$ (resp. $\PP ''$). As follows from a remark of Clemens
on p. 211 in \cite{C-1}, the fiber ${\FFF}_0$ is
acquired in the family $\F /\Delta$ with multiplicity 1.
The lines parametrized by
${\FFF}_0$, that is, the bisecants to $C$ in either one of the
two copies of $\PP^3$, will be called {\em lines in} $X_0$.
The Neron model $\J$ is a family of abelian Lie groups such that
the fiber $J_0$ over $t=0$ is the union of two
components $J^+,J^-$, both isomorphic to a $\CC^*$-extension
of the Jacobian $J(B)$ of $B$. The Neron model has a
natural compactification $\overline{\J}$, obtained by pasting in
two copies of $J(B)$, so that the central fiber
$\overline{J}_0$ is the union of two components $\overline{J}^\pm$
which are both $\PP^1$-bundles over $J(B)$ and which
meet each other transversely along
two disjoint sections of the $\PP^1$-bundles.

In \cite{C-2}, Clemens studies the
family of sextics of genus 3 in
the varieties $X_t$. We will treat in
the same way the one of elliptic quintics.

Let $\bH\lra \Delta$ be the subscheme of the relative Hilbert
scheme $\Hilb (\X /\Delta )$, such that the fiber $H_t$ over $t\neq 0$ is the
component of $\Hilb (X_t)$ whose
generic point represents an elliptic quintic in $X_t$,
and the points of $H_0$, the fiber of $\bH$ over $t=0$,
represent the curves which are limits of elliptic quintics in $X_t$
as $t\rar 0$, that is, $\bH$ is the closure
of $\bigcup_{t\in\Delta^*}H_t$ in $\Hilb (\X /\Delta )$.
Let $\tilde\Delta\lra\Delta$, $s\mapsto t=s^e,$ be some base change and
$C_s\in H_{s^e}$ a family of curves in $X_{s^e}$ which are smooth
elliptic quintics for $s\neq 0$. Their images
$\uC_s=\pi_t(C_s)\subset \PP^3$ are elliptic
quintics which are totally
tangent to $W_t$.
Each $C_s$ for $s\neq 0$
defines a unique K3 surface $S(C_s)\in |\OOO_{X_t}(2)|$
such that $C_s\subset S(C_s)$,
where we denote by $\OOO_{X_t}(k)$ the sheaf $\pi_t^*(\OOO_{\PP^3}(k))$.
The image $\pi_t(S(C_s))$ is a quartic $\underline{S(C_s)}$ in $\PP^3$,
containing $\uC_s$ and
tangent to $W_t$ along the curve $\underline{S(C_s)}\cap W_t$,
so that $\pi_t^{-1}(\underline{S(C_s)})$ is the union of
two components: $S(C_s)$ and $\iota S(C_s)$, where $\iota$
denotes the Galois involution of $\pi_t$.  For $s=0$,
there is a unique limit quartic in $\PP^3$,
possibly singular, which we
denote by  $\underline{S(C_0)}$.
So the family $S(C_s)$ extends over $s=0$
with fiber $S(C_0)\simeq\underline{S(C_0)}$
lying in one of the two components $\PP',\PP''$ of $X_0$,
and the limiting curve $C_0$ lies entirely in $\PP'$ or in $\PP''$.
As remarks Clemens, a component of $C_0$ of degree $d$ which
does not lie in $Q=\PP'\cap\PP ''$ has to meet
$B=W\cap Q$ in $2d$ points counting multiplicities,
because otherwise the deformed curve $\uC_s$
for small $s$ will not be totally tangent to $W_t$, which is absurd.
We will call such components 2d-secants to $B$.

Let $H'$, resp. $H''$ be the subset of $H_0$ consisting of the curves
$C\subset X_0$ such that:
(1) $C$ is smoothable in $X_0$, (2) $C_{\red}$ has no components contained
in $Q$, and (3) all the components of $C$ lie in $\PP '$,
resp. in $\PP ''$. Denote $H^0=H'\cup H''$.

As ${\FFF}_0$ is a reduced
fiber of $\F$, there exists a local analytic family of lines $l_t\in {\FFF}_t$
near $t=0$, and we can use $5l_t$ as the section of
reference curves for the definition of a relative Abel-Jacobi map
$\aalpha =\aalpha_{{\bH}}:\bH\lra\overline{\J}$.
According to \cite{C-1}, Sect. 8, $\aalpha$ extends as
a regular map to all of $H^0\subset H_0$, and the image
of $H'$ (resp. $H''$) lies in
one  component $J'$  (resp. $J''$) of the smooth locus
$J_0=J^+\cup J^-$ of $\overline{J}_0$ (thus, $J'=J^+$
or $J^-$, and $J''$ may a priori coincide either with $J'$, or
with the other component of $J_0$).
If we consider the composition $\phi$ of $\aalpha|_{H^0}$
with the natural projection $\eta :J_0\lra J(B)$,
we obtain the following description
modulo a constant translation by a fixed divisor on $B$ depending
on the choice of the reference point:
$\phi (C)$ is the class of $B\cdot C$ considered as a degree-10
divisor on $B$. In particular, $H'$ dominates $J(B)$
via $\phi$ if and anly if $H''$ does.

We will use these considerations to prove the following proposition.

\begin{proposition}\label{zavyazka}
Let $\ppi :\X\lra\Delta$ be a family of
degenerating quartic double solids as above which is generic, that is,
$F$ and $G$ are generic. Then
there exists a morphism
$$
\aalpha =\aalpha_{\bH} : \bH\setminus A\lra \overline{\J}
$$
over a possibly shrinked disc $\Delta$,
where $A=H_0\setminus H^0$ is a proper closed
subset of $H_0$,
such that the following properties are verified:

(i) The restriction \mbox{$\alpha_t: H_t\lra J_t$} of $\aalpha$
to any fiber $H_t$ with $t\in\Delta^*$ is, modulo a constant
translation in $J_t$, the Abel--Jacobi map of $H_t$.

(ii) The generic fiber of $\alpha_t$ for any $t\in\Delta$
is the union of finitely many copies of $\PP^1$
corresponding to elliptic pencils in the K3 surfaces
$S(C_t)$.

(iii) The closure $\D $ of $\aalpha (\bH \setminus A)$ is a relative
divisor in $\overline{\J}$ over $\Delta$,
that is $D_t=\D\cap \overline{J}_t$ is a divisor in $\overline{J}_t$ for every
$t\in\Delta$, and $D_t$ is a translate of the theta divisor
of $J_t$ if $t\in\Delta^*$.

(iv) The smooth locus $J_0$ of the central fiber $\overline{J}_0$
of $\overline{\J}$ consists of two components $J^\pm$, and only
one of them, say $J'$,
contains $\alpha_0(H^0)$.
Moreover, $\alpha_0(H^0)=\alpha_0(H')=\alpha_0(H'')\subset J'$.

(v) $\alpha_0(H^0)$ is a section
of the natural projection $J'\lra J(B)$ over an open subset \mbox{of $J(B)$}.

\end{proposition}

\begin{proof}
(i) and (ii) follow from the above dixcussion, (iii) is a consequence of
Theorem \ref{Theta}.
To prove (iv) and (v), we will use
the characterization of the divisors in $\overline{J}_0$
which can be limits of
a family of principal polarizations in the neighboring fibers $J_t$,
given by Clemens in \cite{C-2}, (4.2.1)-(4.2.2).
It states that: 1) only one of the irreducible components of
$D_0$, say $D'$, dominates $J(B)$,  and 2) if $\overline{J}^\epsilon$
($\epsilon=+$ or $-$) is the component of $\overline{J}_0$
containing $D'$, then $D'$ is a section
of the natural projection $\overline{J}^\epsilon\lra J(B)$
over an open subset of $J(B)$.

As we have remarked, $H'$ dominates $J(B)$
via $\phi$ if and only if $H''$ does. Hence, if we assume that
$H'$ dominates $J(B)$, then
$\alpha_0(H')$ and $\alpha_0(H'')$ both dominate $J(B)$,
hence $J'=J''$ and $\alpha_0(H')=\alpha_0(H'')$ is
one and the same rational section of $J'\lra J(B)$.

It remains to see that $H'$ dominates $J(B)$. This follows
from the fact that for generic $W,Q$ and for a
generic 10-uple of points $Z$ in $B=W\cap Q$, there is
at least one elliptic quintic in $\PP^3$ passing through $Z$
(see Corollary \ref{10ptsonQ}). Alternatively, one can
prove, by an argument similar  to that of Clemens on p.
97 of \cite{C-2}, that the adherence values of $\aalpha$
along $A$ project to a proper closed subset of $J(B)$,
and thus at least one of the components of $H_0\setminus A =H'\cup H''$
has to project to a rational section of
$\overline{J}^\epsilon\lra J(B)$.
\end{proof}

\end{sub}

\begin{sub}\rm {\bf Deforming pentagons from $X_0$ to $X_t$.} \end{sub}

We want now to prove that the fiber $H_0$ is reduced at
a generic point of $H'$. To this end,
we will construct a local cross-section
of the projection $\bH\lra\Delta$ in the neighborhood of $0$
taking its values in the subschemes of $H_t$ parametrizing
cycles of five lines, or pentagons in $X_t$.

\begin{definition}
Let $X$ be a quartic double solid. A pentagon in $X$ is
a reducible curve $\Gamma =\bigcup_{i=1}^5\ell_i$, where $\ell_i$
are lines in $X$ such that $\ell_i\cap\ell_{i+1}\neq\empt$
(addition of subscripts modulo $5$). If $\ell_i$ meets $\ell_{i+1}$
quasi-transversely at one point
and these are the only intersection points of $\ell_i,\ell_j$
($i\neq j$), we will call $\Gamma$ a good pentagon.
\end{definition}

Consider the 5-uple fibered product of $\F$ with itself
over $\Delta$:
$$
\Y=\F_{\Delta}^5=\F\times_\Delta \F\times_\Delta \F\times_\Delta \F\times_\Delta \F
$$

Let $Y_t$ denote the fiber of $\Y$ over $t\in\Delta$ and
$\D_i$
the subvariety of $Y$ parametrizing the
5-uples $(\ell_1,\ldots, \ell_5)$ of lines in the $X_t$
defined by the incidence relation
$\ell_i\cap\ell_{i+1}\neq\empt$ ($i=1,\ldots, 5$, $5+1=1$).
Then the variety parametrizing all the pentagons
in $X_t$, $t\in\Delta$ can be represented as the
intersection
$\SSigma:=\bigcap_{i=1}^5 \D_i$.

\begin{proposition}\label{pentagons-transv}
Let the quartic $W$ and the quadric $Q$ defining the family
$\X\lra\Delta$ be generic.
Let $\Sigma_0$ be the fiber of $\SSigma$ over $t=0$. Then
there exists a component of $\Sigma_0$, such that
its generic point $y$ represents a good pentagon in $X_0$
whose five vertices are not in $Q$,
the $6$ subvarieties $Y_0$, $\D_1,\ldots ,\D_5$ are nonsingular at $y$
and their intersection is transversal at $y$.
\end{proposition}

\begin{proof}
It suffices to prove the assertion for a special $y_0\in Y_0$ and
a special degenerate quartic $W$. Choose
for $y_0$ the pentagon $\bigcup_{i=1}^5\ell_i^0$, where
$\ell_i^0=\overline{P_iP_{i+1}}$ (addition of subscripts modulo 5),
$P_1=(1:0:0:0)$, $P_2=(0:1:0:0)$, $P_3=(0:0:1:0)$,
$P_4=(0:0:0:1)$, $P_5=(1:1:1:1)$. Let
$W$ be the union of 4 planes $\Pi_k$, $k=1,\ldots , 4$,
where $\Pi_1=\{ x_3-\lambda x_4=0\}$,
$\Pi_2=\{ x_1-\mu x_4=0\}$,
$\Pi_3=\{ x_1-\nu x_2=0\}$,
$\Pi_4=\{ x_1-x_3-\eta (x_2-x_3)=0\}$, and
$Q$ a generic quadric passing through the two
points $A_1^0=\Pi_2\cap \ell_5^0$ and $A_2^0=\Pi_3\cap \ell_5^0$, so that
$C=W\cap Q$ is the union of 4 conics $C_k=\Pi_k\cap Q$, $k=1,\ldots , 4$.
An affine neighborhood of $y_0$ in $\Y$ can be parametrized by
11 rational functions $t,u_1,v_1,\ldots ,u_5,v_5$, where
$(u_i,v_i)$ for $i=1,\ldots ,4$ are affine coordinates
in the plane $\PP^2_i$ parametrizing lines in $\Pi_i$, and $u_5,v_5$
are the parameters representing  $A_1^0\in C_2$ and $A_2^0\in C_3$
respectively. Writing out the linearized equations of $\D_i$
at $y_0$ (depending on $\lambda, \mu, \nu, \eta$ and two
more constants defining tangent directions of $C_2, C_3$
at $A_1^0$, $A_2^0$), one immediately verifies
that they do not depend on $t$ and are linearly independent
as soon as $\mu\neq 0,\nu -\mu\neq 0, \nu -1\neq 0$,
which ends the proof.
\end{proof}
\smallskip

\begin{sub}\rm {\bf Smoothability of good pentagons.}\smallskip

Now we will check that the good pentagons can be deformed
into elliptic quintics.

\begin{proposition}\label{pentagons-smthbl}
Under the hypotheses and in the notation of Proposition \ref{pentagons-transv},
there exists
a local cross-section $\Delta\lra \bH$,
$t\mapsto y(t)$, over a possibly shrunk disc $\Delta$,
such that $y_0\in H'$ and the corresponding curves $\Gamma_{y(t)}\subset X(t)$
are good pentagons, strongly smoothable into elliptic quintics in $X_t$
for all $t\in\Delta$.
\end{proposition}

Recall that $\Gamma_0=\Gamma$ is called {\em smoothable}
if there exists a flat family of curves
$\Gamma_s$ such that its generic member is
smooth. If there exists such a deformation
whose base and total space are both smooth,
then $\Gamma_0$ is {\em strongly smoothable}.
\smallskip

{\em Proof of Proposition} \ref{pentagons-smthbl}.
As above, we specialize
our family $\X\lra\Delta$
to the case when $C$ is the union of four conics $C_i$.
Let $\Gamma =\bigcup_{i=1}^5\ell_i^0$ be a good pentagon as
constructed in the proof of Proposition \ref{pentagons-transv}.

The following lemma is proved by a standard application
of the techniques of \cite{HH}.

\begin{lemma}\label{pentagons-dim20}
The Hilbert scheme of $X_0=\PP '\cup\PP ''$ is smooth
of local dimension $20$ at $\Gamma$.
\end{lemma}

Now we will impose on the curves parametrized by
$\Hilb (X_0)$ the incidence conditions with the four
conics $C_i$.

\begin{lemma}\label{UC}
In the situation of Lemma \ref{pentagons-dim20}, let $U_C$
be the locally closed subset of $\Hilb (X_0)$ whose points
represent the curves meeting transversely $C$
at $10$ points. Then the point $y_0$ representing $\Gamma$
is a smooth point of $U_C$, $\dim_{y_0}U_C=10$,
and the generic point of the component of $U_C$
contaning $y_0$ represents a nonsingular curve in $\PP '$.
\end{lemma}

\begin{proof}
By Lemma \ref{pentagons-dim20}, $y_0$ is a smooth
point of $\Hilb (X_0)$. It has a neighborhood $U$ consisting only
of smooth points and contained
entirely in $\Hilb (\PP ')$. Let $F\subset\PP '$ be any reduced curve
meeting $\Gamma$ transversely at $r$ points $x_1,\ldots ,x_r$ which are
nonsingular both on $F$ and on $\Gamma$, and let $\xi_i$ for
$i=1,\ldots ,r$ denote the point of $\PP (\NNN_\Gamma )$
corresponding to the direction of the branch of $F$ passing through $x_i$.
Denote by $\elm^-_{\xi_1,\ldots ,\xi_r}\NNN_\Gamma$ the negative elementary
transformation of the normal bundle $\NNN_\Gamma$ w.r.t. the $\xi_i$
(see \cite{HH}). Let $U_F\subset U$ be the
locus of curves meeting $F$ transversely at $r$ distinct points
and represented by elements of $U$. Then $U_F$
is smooth at $y_0$ if
$h^1(\Gamma ,\elm^-_{\xi_1,\ldots ,\xi_r}\NNN_\Gamma )=0$
and in this case we have for the tangent space $T_{y_0}U_F=H^0(\Gamma ,
\elm^-_{\xi_1,\ldots ,\xi_r}\NNN_\Gamma )$
(see \cite{Ran}, Remark 5.1).

Let us apply this observation to $F=C=\bigcup_{i=1}^4C_i$, $r=10$,
$C_i\cap \ell_i= \{ x_{2i-1},x_{2i}\}$, $i=1,\ldots ,5$.
Let $\NNN = \elm^-_{\xi_1,\ldots ,\xi_{10}}\NNN_\Gamma $. We have
$\NNN_{\Gamma}|_{\ell_i}\simeq 2\OOO (2)$ for any $i$.
As the branches of $C$ passing through $x_{2i-1},x_{2i}$
lie in one plane for each $i=1,\ldots ,4$ and they do not
for $i=5$, we have $\NNN |_{\ell_i}\simeq \OOO\oplus \OOO (2)$
for $i=1,\ldots ,4$ and $\NNN |_{\ell_5}\simeq 2\OOO (1)$.
This easily implies that $h^0(\NNN )=10$, $h^1(\NNN )=0$.
It remains to prove that $\Gamma$ is smoothed out by
a generic deformation inside $U_C$. It suffices to
show that the natural map $T_{y_0}U_C\lra T^1_{p_i}$ to the Schlesinger's
space $T^1$ is surjective for any $i=1,\ldots ,5$, where
$p_i=\ell_i\cap\ell_{i+1}$ (the addition of subscripts modulo 5). This
can be done
in the same way as  in the proof
of Theorem 4.1 in \cite{HH}.
\end{proof}

Now we go back to the proof of Proposition \ref{pentagons-smthbl}.
By Proposition \ref{pentagons-transv}, we can find a local
cross-section of $\Sigma$ representing a flat
family of good pentagons in the varieties $X_t$, whose fiber
over $0$ is a good pentagon in $\PP '$. By the universal
property of the Hilbert scheme, this family is induced via
a cross-section $\Delta\lra\Hilb (\X /\Delta )$ over $\Delta$,
which we will denote by $t\mapsto y(t)$. It suffices to show that
$y_0\in H'$. We have verified in Lemma \ref{UC}
that a good pentagon $\Gamma$ in $\PP '$ becomes
strongly smoothable when $C$ is specialized into a union of four
conics. By the semi-continuity of $h^i(\NNN )$, and because
the surjectivity of the maps $H^0(\NNN )\lra T^1_{p_i}$
is an open condition, a generic good pentagon
in some component of good pentagons in $\PP '$ is
strongly smoothable for generic $\X\lra \Delta$.
In fact, we do not check that
the family of good pentagons in
$\PP '$ is irreducible, so we
are speaking here about one of its components.

As $\Gamma$ is a connected curve of degree 5 and arithmetic
genus 1, so is its smoothing, hence it is an elliptic quintic in
$\PP '$. By construction, it meets $C$ transversely at 10 points
and is an element of $H'$.
This ends the proof of Proposition \ref{pentagons-smthbl}.
\square

This immediately implies

\begin{corollary}\label{Hmult1}
With the same assumptions as above, there is a component of
$H'$ which is acquired with multiplicity $1$ in the fiber $H_0$
of $\bH\lra \Delta$.
\end{corollary}

We will prove in Section \ref{sec-etale} that $H'$ is
irreducible, hence $H'$ and $H''$ are irreducible
components of $H_0$ acquired with multiplicity 1.
\end{sub}

\section{Elliptic quintics through ten points}
\label{sec-etale}

\begin{sub}\rm {\bf
Elliptic quintics, 10-secant to $\boldsymbol B\boldsymbol =\boldsymbol W
\boldsymbol \cap \boldsymbol Q$.}

We are going to prove the irreducibility
of the subscheme $H'$ of $\Hilb (X_0)$ introduced in Sect. \ref{deg-pairs}.

\begin{lemma}\label{H-prime-irred}
Let $Q$ (resp. $S$) be a generic quadric (resp. quartic) in
$\PP^3$, and $B=Q\cap S$. Let $\CCC_5^1[B]_{10}\subset \CCC_5^1(\PP^3)$ be
the base of the family of elliptic quintics meeting
$B$ in a subscheme of length $10$. Then $\CCC_5^1[B]_{10}$ is irreducible.
\end{lemma}

\begin{proof}

Let $\PP^{34}$ be the space of quartic surfaces in $\PP^3$,
$\Gamma=\{(S,C)\in\CCC^1_5(\PP^3)\times
\PP^{34}~|~C\subset S\}$ and $\PP^{34}\xleftarrow{p_1}\Gamma\xrightarrow{p_2}
\CCC^1_5(\PP^3)$ the natural projections.
Any $C\in\CCC^1_5(\PP^3)$ does not lie on a
quadric, hence the restriction map $H^0(\OOO_{\PP^3}(2))\to H^0(\OOO_C(2))$
is bijective. Then
$h^1(\III_{C/\PP^3}(2))=0$ and a fortiori
 $h^1(\III_{C/\PP^3}(4))=0$, which implies $h^0(\III_{C/\PP^3}(4))=15.$
Hence $p_2^{-1}(C)\simeq\PP^{14}$, and as
$\CCC^1_5(\PP^3)$ is irreducible of dimension 20,
$\Gamma$ is irreducible as well
and $\dim\Gamma=34$. Moreover, denoting
$\CCC^1_5(\PP^3)^*:=\CCC^1_5(\PP^3)\smallsetminus\Sing(\CCC^1_5(\PP^3)),\ \
\Gamma^*=p_2^{-1}(\CCC^1_5(\PP^3)^*)$,
we see that $\Gamma^*$ is smooth and is a
projective bundle over $\CCC^1_5(\PP^3)^*$.

Next, a general quartic $S$ in
$p_2^{-1}(C)$ is smooth for general
$C\in\CCC^1_5(\PP^3)$.
Hence $D=p_1(\Gamma)$ is
irreducible. Besides, for any
smooth quartic $S=p_1(S,C)\in D$ the fibre
$p_1^{-1}(S)$ is a dense open
subset in $|\OOO_S(C)|\simeq\PP^1$;
hence $\dim D=33$, i.e. $D$ is an
irreducible divisor in $\PP^{34}$.
Moreover, $p_1:\Gamma\to D$ is a smooth morphism of relative dimension 1.

Let $\overline{ D}$ be the closure of $D$ in $\PP^{34}$. We will prove that
$
\codim_{\overline D}\Sing\overline D\ge2.
$
In fact, if $\codim_{\overline D}\Sing\overline D=1,$ i.e.
$\dim\Sing\overline D=32$, then
$p_1(p_2^{-1}(C))\cap\Sing\overline D=\PP^{14}\cap
\Sing\overline D\ne\emptyset$ for any
$C\in\CCC^1_5(\PP^3)^*$, hence, by the smoothness of
$p_1|\Gamma^*$, $p_2^{-1}(C)\cap
\Sing\Gamma^*\subset p_2^{-1}(C)\cap p_1^{-1}
(\Sing\overline D)\ne\emptyset$, which contradicts the
smoothness of $\Gamma^*$.

Now let $\CCC_{2\cdot 4}(\PP^3):=\{C\in\Hilb(X)~|~C$
be a smooth complete intersection
curve of a smooth quartic surface and a smooth quadric
surface in $\PP^3\}$. Clearly,
$\CCC_{2\cdot 4}(\PP^3)$ is a smooth irreducible variety
of dimension 33. Let $U$ be the
open subset of smooth quartics in $\PP^{34}$,
$\Pi=\{(C,F)\in\CCC_{2\cdot 4}(\PP^3)\times\PP^{34}~|~C\subset F\}$
and
$\CCC_{2\cdot 4}(\PP^3)\xleftarrow{q_1}\Pi\xrightarrow{q_2}\PP^{34}$
the projections.
One immediately verifies, using the resolution $0\to\OOO_{\PP^3}(-2)\to
\OOO_{\PP^3}(2)\oplus\OOO_{\PP^3}\to
\III_{C/\PP^3}(4)\to0$ for $C\in\CCC_{2\cdot 4}(\PP^3)$,
that $h^0(\III_{C/\PP^3}(4))=11$.
By the definition of $\CCC_{2\cdot 4}(\PP^3)$,
a general quartic in
$q_1^{-1}(C)$ is smooth, hence $q_1^{-1}(C)$ is an
open subset in $\PP^{10}=P(H^0(\III_{C/\PP^3}(4))).$ Thus,
$\Pi$ is irreducible of dimension 43.
Moreover, for any $(C,S)\in\Pi$, we have $(C,S)\in q_2^{-1}(S')$
with $S'\in U$, and
$q_2^{-1}(S')$ is an open subset in $|\OOO_{S'}(C)|=|\OOO_{S'}(2)|\simeq
\PP^9$ consisting of smooth curves. Clearly, $D':=U\cap D\ne\emptyset$ and
$D'$ is open
in $D$, hence $D_\Pi:=q_2^{-1}(D')$ is an
irreducible divisor in $\Pi$.

Let us look now at the general fiber of
$q_D:=q_1|D_\Pi:D_\Pi\to\CCC_{2\cdot 4}(\PP^3)$.
Let $\overline\Pi$ be the closure of $\Pi$ in
$\overline\CCC\times\PP^{34}$ with natural
projections $\CCC\xleftarrow{\overline q_1}\overline\Pi
\xrightarrow{\overline q_2}\PP^{34}$, where
$\overline\CCC$ is the closure of $\CCC_{2\cdot 4}(\PP^3)$
in $\Hilb(\PP^3)$, and  $\overline D_\Pi =\overline q_2^{-1}(\overline D)$.
Let $Z\subset\overline D$ be the subset of $\overline D$
which consists of reducible quartics,
$\codim_{\overline D}Z\ge11$, $D^*=\overline D\smallsetminus Z$,
$D^*_\Pi:=\overline q_2^{-1}(D^*)$.
For any
$S\in D^*$, $|\OOO_S(2)|\simeq\PP^9$ lies in
$\overline\CCC$. Hence
$\overline q_2|_{D^*_\Pi}:D^*_\Pi\to D^*$
is a projective bundle with fiber $\PP^9$,
in particular, it is a smooth morphism of relative
dimension 9. By construction, the fiber of
$\overline q_D=\overline q_1|\overline D:
\overline D_\Pi\to\overline\CCC$ over any
$C\in\CCC_{2\cdot 4}(\PP^3)$  is a divisor in
$\overline q_1^{-1}(C)=\PP^{10}$ which contains
$q_D^{-1}(C)$ as an open subset.
Assume that $q_D^{-1}(C)$ is reducible. Then its
closure $\overline{ q_D^{-1}(C)}$ is also a reducible
divisor in $\PP^{10}$, hence it has
singularities in codimension 1. Since $C\in\CCC_{2\cdot 4}(\PP^3)$
is smooth, $\codim_{\overline D_\Pi}\Sing\overline D_\Pi=1$.
As $\codim_{\overline D}Z\ge11$,  $\codim_{ D^*_\Pi}\Sing
 D^*_\Pi=1$. Hence, since
$\overline q_2|_{ D^*_\Pi}:D^*_\Pi\to D^*$ is a smooth morphism,
$\codim_{ D^*}\Sing D^*=1$ and a
fortiori $\codim_{\overline D}\Sing\overline D\ge2$,
which is absurd.

Thus, for general $C\in\CCC_{2\cdot 4}(\PP^3)$,
both $q_D^{-1}(C)$ and $\overline{ q_D^{-1}(C)}$ are irreducible.
The generic quartic $S$ from $\overline{ q_D^{-1}(C)}$ is smooth
and contains a unique pencil $|E|$
of elliptic curves, and these curves
belong to $\CCC_5^1[B]_{10}$, where $B=Q\cap S$, $Q$ being the unique
quadric containing $C$.
Conversely, taking $B$ as in the hypothesis of the lemma,
any smooth elliptic quintic $E$
meeting $B$ in $10$ points lies in a
generically smooth quartic surface $S$
through $B$, i.e. such that $(B,S)\in\overline q_D^{-1}(B)$,
and varies on $S$ in a pencil.
Hence we have a fibration $H'\to\overline q_D^{-1}(C)$ with an open subset of $\PP^1$ as a
fiber. This implies the irreducibility of $\CCC_5^1[B]_{10}$.
\end{proof}

\begin{corollary}\label{Hirred}
Let $Q\subset\PP^3$ be a smooth quadric surface,
$W\subset\PP^3$ a generic quartic, $B=Q\cap W$.
Then the subscheme $H'$ (resp. $H''$) of $\Hilb (X_0)$ introduced in Sect. \ref{deg-pairs}
and $\CCC_5^1[B]_{10}$ contain isomorphic dense open subsets, hence
$H'$ (resp. $H''$) is irreducible. This implies, by Corollary \ref{Hmult1},
that the fiber $H_0$ of $\bH$ over $t=0$ is reduced at
the generic points of $H'$ and $H''$.
\end{corollary}

\end{sub}

\begin{sub}\rm   {\bf 42 elliptic quintics through 10 points}\label{42}

Let $U$ be the open subset of the smooth locus of $\Hilb\PP^3$
parametrizing elliptic quintics.
There is a universal family
$\CCC\subset U\times\PP^3$ over $U$ with
projection $p_\CCC:\CCC\to U$.
Consider the 10-th relative symmetric power $\SSS:=\Sym^{10}(\CCC/U)$ of $\CCC$ over $U$
with projection $p:\SSS\to U$ and fiber
$
p^{-1}(\{C\})=\Sym^{10}C.
$
Thus $\SSS$ is a smooth irreducible variety of dimension 30.
We can embed it into $U\times \Hilb^{10}\PP^3$:
$$\SSS=\{(\{C\},\{Z\})~|~\{C\}\in U,\{Z\}\in\Hilb^{10}\PP^3
\ \mbox{and}\ Z\subset C\ \mbox{as a scheme}\} ,$$
and we have a natural projection
$
\mu:\SSS\to\Hilb^{10}\PP^3,\ (\{C\},\{Z\})\mapsto\{Z\}.
$
There exists at least one
quintic from the closure of $\CCC_5^1(\PP^3)$ passing
through 10 distinct points in $\PP^3$, and for a
general 10-uple of points, the number of
such quintics is finite and they are smooth.
In fact, we have:
\begin{lemma}[E.~Getzler, \cite{Ge}, Table 2]\label{ezra}
For a general 10-uple $Z$ of points in $\PP^3$,
there are $42$ elliptic quintics passing through $Z$.
\end{lemma}

Denote by $H^{10}$ the component of $\Hilb^{10}\PP^3$
containing a 10-uple of distinct points.
It follows that $H:=\mu(\SSS)$ is an open subset of $H^{10}$ and that
$\mu:\SSS\to H$ is generically finite of degree 42.
We need to replace in the above lemma
``a general 10-uple $Z$ of points in $\PP^3$"
by ``a general 10-uple $Z$ of points in a smooth quadric $Q\subset\PP^3$",
that is, we should prove that the 10-uples of points lying in
smooth quadrics form a divisor in $H$ which is not contained in
the branch locus of $\mu$.

Let
$D=\{(\{C\},\{Z\})\in\SSS~|~Z=C\cap Q\
\mbox{for some quadric $Q$ in}\ \PP^3\} .$
Clearly, $D$ is a divisor in $\SSS$, and for any
$\{C\}\in U$,
\begin{equation}\label{DcapFiber}
D\cap p^{-1}(\{C\})=\{\{Z\}\in \Sym^{10}C~|~\OOO_C(Z)=
\OOO_{\PP^3}(2)|C\}=\{\big|\OOO_{\PP^3}(2)|C\big|\}\simeq\PP^9.
\end{equation}
As $U$ is irreducible,
$D$ is irreducible by (\ref{DcapFiber}), hence so is
$B:=\mu(D)$. It is easy to see that $B$ is a divisor.
Indeed, for any smooth elliptic quintic $C$ in $ \PP^3$
and for any $Z\in\:|\OOO_C(2)|$, there is a unique quadric $Q$
passing through $Z$, so $\dim B=29$.

\begin{proposition}\label{etale}
Let $R\subset \SSS$ be the ramification divisor of $\mu$.
Then we have:
(i) $D$ is not contained in $R$, hence
$\mu: \SSS\to H$ is \'etale in a general point of $D$.

(ii) $\mu$ is an \'etale covering over an open subset $V$ of $H$ such that
$B\cap V\ne\empt.$
\end{proposition}

\begin{proof} First we show that (i) implies (ii).
For a general 10-uple of points
$Z$ on a general smooth quadric $Q$, there exists a quartic surface
passing through $Z$ and intersecting $Q$ in a
smooth curve, say $C_0$. Further, by Corollary \ref{Hirred},
the set of all the elliptic quintics meeting $C_0$
in 10-uples of points is irrreducible. Hence
$(\mu^{-1}(\Sym^{10}(C_0)))_{\red}$
is irreducible. Hence $\mu$ has the same ramification indices in
all the points of
$\mu^{-1}({Z})$, where $\{Z\}\in B$ is the above general point of $B$.
But according to
(i), $\mu$ is
unramified in at least one point of $\mu^{-1}({Z})$.
Hence $\mu$ is unramified, that is \'etale, in all
the points of $\mu^{-1}({Z})$.

Let us prove (i).
A pair $(Z,C)$ is a ramification point of $\mu$
if and only if $C$ has infinitesimal deformations
which fix $Z$, that is, $h^0(\NNN_{C/\PP^3}\otimes\III_Z)\neq 0$.
But if $Z$ lies on a quadric, $\NNN_{C/\PP^3}\otimes\III_Z
\simeq \NNN_{C/\PP^3}(-2)$, and $h^0(\NNN_{C/\PP^3}(-2))=0$
by \cite{EL} (see also Theorem (VIII.2.7) in \cite{Hu}).
\end{proof}

The following statements are obvious corollaries of what we
have proved:

\begin{corollary}\label{10ptsonQ}
Let $Q\subset\PP^3$ be a smooth quadric surface.
Then there are $42$ elliptic quintics in $\PP^3$
passing through $10$ generic points of $Q$.
There are also $42$ elliptic quintics in $\PP^3$
passing through $10$ generic points on the curve $B=Q\cap W$
for a generic quartic surface $W\subset \PP^3$.
\end{corollary}


\begin{corollary}\label{84-0}
Let $\alpha_0:H'\cup H''\lra J^+$
be the map introduced in Proposition \ref{zavyazka}, (ii).
Then the number of copies of $\PP^1$ in
the generic fiber of $\alpha_0$ is equal to $84$.
\end{corollary}

As $H'\cup H''$ contains all the elliptic quintics which are
the limits of elliptic quintics from the neighboring fibers,
we deduce from Corollaries \ref{84-0} and \ref{Hirred}
that the number of $\PP^1$'s in the generic fiber
of $\alpha_t$ is 84 for all sufficiently
small $t\in\Delta$. This implies the following result.

\begin{theorem}\label{84-X}
Let $X$ be a general quartic double solid.
Then the generic fiber of the Abel Jacobi map
$\Phi_{\CCC_5^1(X)}:\CCC_5^1(X)\lra \Theta+\const$ is a dense open
subset in the union of
$84$ copies of $\PP^1$.
\end{theorem}

\end{sub}

\section{Vector bundles}
\label{vb}

Let $X$ be general, so that $\CCC^1_5(X)$ is irreducible. Let
$C\in  \CCC^1_5(X)$. The Serre construction allows us to define an
$\OOO_X$-sheaf
$\EEE_C$
as an extension of $\OOO_X$-sheaves:
\begin{equation}\label{extn}
0\to\OOO_X(-1)\to\EEE_C\to\JJJ_{C/X}(1)\to0.
\end{equation}

Any such extension is determined uniquely up to an isomorphism by a
1-dimensional subspace of the group
$\Ext^1(\JJJ_{C/X}(1),\OOO_X(-1))\simeq H^0(\EXT^1(\JJJ_{C/X}(1),\OOO_X(-1)))
\simeq H^0(\EXT^2(\OOO_C(1),\OOO_X(-1)))\simeq H^0(\EXT^2(\OOO_C,\omega_X))
\simeq H^0(\omega_C)\simeq\CC$
(we are using
$h^i\OOO_X(-2)=0,\ \ i<3,$
and $\omega_C\simeq\OOO_C$).
Thus $\EEE_C$ is determined uniquely up to an isomorphism by $C$.
Besides, by a well known property of the Serre construction,
the local triviality of  $\EXT^1(\JJJ_{C/X}(1),\OOO_X(-1))\simeq
\omega_C$ implies that $\EEE_C$ is locally free.

From (\ref{extn}), one easily computes the Chern classes:
$$
c_1(\EEE_C)=0,\qquad c_2(\EEE_C)=3
$$
(here we understand $c_2$ as an element of the group $B_1$
of classes of algebraic equivalence of 1-dimensional cycles on $X$,
generated by the class of a line on $X$).
As
$\pi:C\to\pi(C)$ is an isomorphism onto a quintic $\pi(C)$
that does not lie in a plane and
$H^0(\OOO_X(1))=\pi^*H^0(\OOO_{\PP^3}(1))$,
it follows that
$H^0(\JJJ_{C/X}(1))\simeq H^0(\JJJ_{\pi(C)/\PP^3}(1))=0$.
Hence the triple (\ref{extn}) gives
$h^0(\EEE_C)=0$.
This means that $\EEE_C$ is Gieseker stable
(see \cite[Prop. 2.6]{MT-1}). Thus, denoting
by $M_X(2;0,3)$ the Gieseker-Maruyama moduli scheme of stable rank-2 vector
bundles on $X$ with Chern classes $c_1=0,\ c_2=3$,
we have $[\EEE_C]\in M(2;0,3)$,
where $[\EEE]$ stands for the isomorphism class of $\EEE$.
Remark that the equality
$h^0(\JJJ_{C/X}(2))=1$
immediately implies
$h^i(\JJJ_{C/X}(2))=0,\ \ i>0;$
hence by (\ref{extn}) also
$$
h^i(\EEE_C)=0,\ \ i>0.
$$
Let
$\tilde M=\{[\EEE]\in M_X(2;0,3)|\ h^i(\EEE)=0,\ \ i>0\}$.
By semicontinuity, this is an open subset of
$M_X(2;0,3)$ containing $[\EEE_C]$. Let
$M$ be the irreducible component of $\tilde M$
containing $[\EEE_C]$. We have a natural map
$$
f: \CCC^1_5(X)\to M,\ C'\mapsto[\EEE_{C'}].
$$
It is a standard matter to show that $f$ is a dominant morphism
(see \cite[Lemmas 5.2 and 5.3]{MT-1}),
whose fibre
$f^{-1}(f(C'))$
is open in
$P(H^0(\EEE_{C'}(1))\simeq\PP^1$:
$$
\overline{f^{-1}(f(C'))}=\PP^1,\ \ \ C'\in  \CCC^1_5(X).
$$
The closure is taken in $\overline{\CCC^1_5(X)}$, and we keep in mind that,
by Riemann-Roch,
$h^0(\EEE_{C'}(1))=2$,
since by the definition of $\tilde M$,
$h^i(\EEE_{C'}(1))=0,\ \ i>0$.
The equality
$h^0(\JJJ_{C'/X}(2))=1$
means that the linear series
$|\JJJ_{C'/X}(2))|$
consists of a unique K3 surface $S(C')$,
and by Serre construction, the fibre
$\overline{f^{-1}(f(C'))}=\PP^1$
is just the pencil of elliptic curves on $S(C')$,
which are nothing else than the zero loci of sections
of $\EEE_{C'}(1)$ :
\begin{equation}\label{pencil}
\overline{f^{-1}(f(C'))}=|\OOO_{S(C')}(C')|=\PP^1.
\end{equation}

Now consider the Stein factorization of the Abel-Jacobi map
$\Phi=\Phi_{\CCC^1_5(X)}: \CCC^1_5(X)\to\Theta+{\rm const.}$
Using the above description (\ref{pencil}) of fibres of $f$
and the well-known interpretation of the differential of the Abel-Jacobi
map in terms of the branch quartic $W$ of
$\pi:X\to\PP^3$ (see the end of Section \ref{abjac}, or
\cite[Corollary 1]{T-1}, or \cite[Prop.2.13]{We}),
we obtain the following

\begin{theorem}\label{Stein}
Let $X$ be a general quartic double solid. Then there exists a
quasi-finite dominant morphism $g:M_0\to\Theta+{\rm const}$
of degree $84$ from
an open part $M_0$ of the above component $M$ of the moduli space
$M_X(2;0,3)$ to
a translate of the theta divisor of $J(X)$
which gives the Stein factorization of the Abel-Jacobi map
$\Phi$ of $ \CCC^1_5(X)$:
\begin{equation}\label{triangl}
\xymatrix{ \CCC^1_5(X)
 \ar[dr]_{f} \ar[rr]^{\Phi} & &
\Theta+{\rm const}\\
 & M_0 \ar[ur]_{g}  & \\ }
\end{equation}
The fibers of
$f: \CCC^1_5(X)\to M_0$
are open subsets of $\PP^1$.
\end{theorem}

The proof is similar to that of Theorem 5.6
in \cite{MT-1}.

\end{document}